\documentclass[reqno]{amsart}
\usepackage{amsmath,amsthm,amscd,amssymb,amsfonts, amsbsy}

\usepackage[usenames, dvipsnames]{color}
\usepackage{enumitem}

\usepackage{mathrsfs}

\usepackage{marginnote}
\usepackage{todonotes}
\usepackage{hyperref}
\usepackage{esint}
\usepackage{verbatim}

\theoremstyle{plain}
\newtheorem{theorem}{Theorem}[section]
\newtheorem{lemma}[theorem]{Lemma}
\newtheorem{corollary}[theorem]{Corollary}
\newtheorem{proposition}[theorem]{Proposition}

\theoremstyle{definition}
\newtheorem{definition}[theorem]{Definition}

\theoremstyle{remark}
\newtheorem{remark}[theorem]{Remark}

\newcommand{\M}{{\mathcal M}}

\numberwithin{equation}{section}

\newcommand{\bM}{\mathfrak{M}}

\newcommand{\bR}{\mathbb{R}}
\newcommand{\bH}{\mathbb{H}}

\newcommand\cD{\mathcal{D}}

\newcommand\cL{\mathcal{L}}

\newcommand\sL{\mathscr{L}}

\newcommand\sW{\mathscr{W}}

\newcommand\sH{\mathscr{H}}

\providecommand{\norm}[1]{\lVert#1\rVert}

\makeatletter
\@namedef{subjclassname@2020}{%
  \textup{2020} Mathematics Subject Classification}
\makeatother

\begin{document}

\title[Singular elliptic and parabolic equations in non-divergence form]{Weighted mixed-norm $L_p$ estimates for equations in non-divergence form with singular coefficients: the Dirichlet problem}

\author[H. Dong]{Hongjie Dong}
\address[H. Dong]{Division of Applied Mathematics, Brown University,
182 George Street, Providence, RI 02912, USA}
\email{Hongjie\_Dong@brown.edu}
\thanks{H. Dong was partially supported by the Simons Foundation, grant no. 709545, a Simons fellowship, grant no. 007638, and the NSF under agreement DMS-2055244.}

\author[T. Phan]{Tuoc Phan}
\address[T. Phan]{Department of Mathematics, University of Tennessee, 227 Ayres Hall,
1403 Circle Drive, Knoxville, TN 37996, USA }
\email{phan@math.utk.edu}
\thanks{T. Phan was partially supported by the Simons Foundation, grant \# 354889.}

\begin{abstract}
We study a class of non-divergence form elliptic and parabolic equations with singular first-order coefficients in an upper half space with the homogeneous Dirichlet boundary condition. In the simplest setting, the operators in the equations under consideration appear in the study of fractional heat and fractional Laplace equations. Intrinsic weighted Sobolev spaces are found in which  the existence and uniqueness of strong solutions are proved under certain smallness conditions on the weighted mean oscillations of  the coefficients in small parabolic cylinders. Our results are new even when the coefficients are constants and they cover the case where the weights may not be in the $A_p$-Muckenhoupt class.
\end{abstract}

\maketitle

%\today
\section{Introduction}

Denote $\Omega_T = (-\infty, T) \times \bR^{d}_{+}$, where $T  \in (0, \infty]$ is a given number, and $\bR^d_+ = \bR^{d-1} \times \bR_+$ is the upper half space with $\bR_+ = (0, \infty)$. For  a point $x \in \bR^d_+$, we write $x = (x', x_d) \in \bR^{d-1} \times \bR_+$. In this paper, we prove the following theorem regarding elliptic and parabolic equations with singular first-order coefficients, in which $L_p(\cD, \omega)$ denotes the weighted Lebesgue space with a given weight $\omega$ in a domain $\cD$, and $D_d$, $D_{x'}$ denote the partial derivatives in the $x_d$-variable and the $x'$-variable, respectively.

\begin{theorem}
                                \label{thm0}
Let $\alpha\in (-\infty,1)$, $p\in (1,\infty), \gamma \in (\alpha p -1, 2p-1)$, and $\lambda>0$.

\textup{(i)} For any $f \in L_p(\bR^d_+, x_d^\gamma dx)$, there exists a unique strong solution $u=u(x)$ to the equation
\begin{equation}  \label{eq1.52}
\left\{
\begin{array}{cccl}
\Delta u+\frac \alpha {x_d} D_d u-\lambda u & = & f &\quad\text{in}\quad\bR^d_+, \\
u & = & 0 & \quad \text{on} \quad \partial \bR^d_+,
\end{array} \right.
\end{equation}
which satisfies
\begin{equation}
                    \label{eq1.11}
\begin{aligned}
&\int_{\bR^d_+}\Big(|DD_{x'}u|^p+ \big|D^2_d u+\alpha x_d^{-1} D_du\big|^p  + |\sqrt{\lambda} Du|^p
 +|x_d^{-1}D_{x'}u|^p \\
&\quad + |\lambda u|^p +|\sqrt\lambda x_d^{-1} u|^p\Big)x_d^{\gamma}\,dx \leq  N \int_{\bR^d_+}|f|^px_d^{\gamma}\,dx,
\end{aligned}
\end{equation}
where $N=N(d,\alpha, \gamma, p)>0$.

\textup{(ii)} For any $f \in L_p(\Omega_T, x_d^\gamma dxdt)$, there exists a unique strong solution $u=u(t,x)$ to the equation
\begin{equation}
                        \label{eq5.51}
                        \left\{
                        \begin{array}{cccl}
u_t-\Delta u-\frac \alpha {x_d} D_d u+\lambda u & = & f & \quad\text{in}\quad \Omega_T, \\
u & = & 0 & \quad \text{on} \quad (-\infty, T) \times \partial \bR^d_+,
\end{array} \right.
\end{equation}
which satisfies
\begin{equation}
                    \label{eq1.11b}
\begin{aligned}
&\int_{\Omega_T}\Big(|u_t|^p+|DD_{x'}u|^p+ \big|D^2_d u+\alpha x_d^{-1} D_du\big|^p+ |\sqrt{\lambda} Du|^p  +|x_d^{-1}D_{x'}u|^p \\
&\quad + |\lambda u|^p +|\sqrt\lambda x_d^{-1} u|^p\Big)x_d^{\gamma}\,dx dt  \leq  N \int_{\Omega_T}|f|^px_d^{\gamma}\,dx dt,
\end{aligned}
\end{equation}
where $N=N(d,\alpha,\gamma, p)>0$.
\end{theorem}

Theorem \ref{thm0} is a special case of Theorems \ref{para-main.theorem} and  \ref{elli-main.theorem}  below, in which more general equations with measurable coefficients and estimates in weighted mixed-norm Sobolev spaces with $A_p(\mu_2)$-Muckenhoupt weights are considered, where $\mu_2(s) =s^{\gamma_0}, s >0,$ with $\gamma_0 \in (-1, 1-\alpha]$. We refer the reader to Section \ref{sec2} for the definitions of function spaces and strong solutions. A novelty of Theorem \ref{thm0} is that  our weight $x_d^{\gamma}$ is in general not an $A_p(\bR, dx)$-Muckenhoupt weight as usually required in the theory of weighted estimates. Observe  that we do not have the control of the $L_p$-norm of $D_d^2u$ as usual, but instead we  control the $L_p$-norm of $D_d^2 u +\alpha x_d^{-1}D_du$. See an intuitive reason for this fact in the paragraph containing \eqref{ODE.exple} below, and also Lemma \ref{s-W-2-p.constant.eqn} and Remark \ref{rm-47}. When $\alpha=\gamma=0$, the estimates \eqref{eq1.11} and \eqref{eq1.11b} are the classical Calder\'{o}n-Zygmund estimates for the Laplace and heat equations in the half space. When $\alpha =0$, weighted estimates similar to these in Theorem \ref{thm0}  were first obtained in \cite{K99}, and the necessity of such results in stochastic partial differential equations is explained in \cite{Krylov-94}. See also \cite{Kr99, KN} and \cite{KKL22} for further results and recent developments on weighted estimates for equations and systems with bounded, measurable, and uniformly elliptic coefficients.
To the best of our knowledge, Theorem \ref{thm0} is new when $\alpha\neq 0$. It is worth noting that the Dirichlet boundary condition is an effective boundary condition only when $\alpha<1$. For example, when $d=\alpha=1$, the equations \eqref{eq1.52} is equivalent to a 2D Poisson equation in the punctuated plane $\bR^2\setminus \{0\}$ with the zero boundary condition prescribed at the origin. It is well known that such boundary condition is negligible as the Brownian motion in 2D is null recurrent.

Elliptic and parabolic equations with singular coefficients emerge naturally in both pure and applied problems. We refer the reader to \cite{Dong-Phan} for some references of related problems in probability, geometric PDEs, porous media, mathematical finance, and mathematical biology. The equations considered in Theorem \ref{thm0} are also closely related to extension operators of the fractional heat and fractional Laplace equations studied, for instance, in \cite{Caffa-Sil, ST17}. In the literature, much attention has been paid to regularity theory for such equations with singular (or degenerate) coefficients. See, for examples, the book \cite{OR} and \cite{Grushin, Gr-Sa, Fabes-1, Fabes, Smith-Stredulinsky} for classical results, and \cite{Lin, WWYZ, Pop-1} for some recent results.  In \cite{MNS21}, the authors obtained interesting $L_p$ type estimates for extension operators with constant coefficients and the Dirichlet or Neumann boundary conditions, by using an functional analytic approach. See also \cite{MNS22} for results about more general operators in the form
$$
x_d^{\alpha_1}\Delta_{x'}+x_d^{\alpha_2}\big(D_{d}^2+cx_d^{-1}D_{d}-bx_d^{-2} ),\quad\alpha_1,\alpha_2\in \bR,
$$
where $b$ and $c$ are constants. We also mention the recent interesting work \cite{Sire-1,Sire-2}, in which the authors obtained H\"{o}lder and Schauder type estimates for scalar elliptic equations of a similar type under the conditions that the coefficient matrix is symmetric, sufficiently smooth, and the boundary is invariant with respect to the leading coefficients that is essentially the same as \eqref{eq8.47} below, even though we do not assume the symmetry condition on the coefficients.

In \cite{Dong-Phan-RMI}, we obtained the Sobolev type estimates for non-divergence form elliptic and parabolic equations similar to \eqref{eq1.52} and \eqref{eq5.51} in a half space with the Neumann boundary condition when $\alpha\in (-1,1)$. The results were later extended in \cite{DP} to more general $\alpha\in (-1,\infty)$, which is optimal.  The corresponding singular-degenerate equations in divergence form were studied in \cite{Dong-Phan-1, DP} with the conormal boundary condition and in \cite{Dong-Phan} with the Dirichlet boundary condition. In these papers, we dealt with leading coefficients which are measurable in the normal space direction and have small mean oscillations in small cylinders (or balls) in time and the remaining space directions. This is called the partially VMO condition and was first introduced in \cite{Kim-Krylov,Kim-Krylov-1} for non-degenerate equations with bounded coefficients.  For non-divergence form equations, it was assumed that
\begin{equation}
                                    \label{eq8.47}
a_{dj}\equiv 0\quad\text{or}\quad a_{dj}/a_{dd}\text{ are constant for }j=1, 2,\ldots,d-1.
\end{equation}
See  \cite[(1.8)]{Dong-Phan-RMI} and \cite[(1.7)]{DP}. We also refer to a related work \cite{MP} in which a conormal boundary value problem for equations  in divergence form with singular-degenerate coefficients as $A_2$-Muckenhoupt weights was considered.

To give a formal description of our main results for general equations, we introduce some notation. Assume that $a=(a_{ij}) :\Omega_T \rightarrow \bR^{d\times d}$ is a matrix of measurable functions that satisfies the following uniform ellipticity and boundedness conditions with the ellipticity constant $\nu>0$:
\begin{equation} \label{ellipticity}
\nu |\xi|^2 \leq a_{ij}(t, x)\xi_i \xi_j \quad \text{and} \quad |a_{ij}(t, x)| \leq \nu^{-1}
\end{equation}
for any $\xi = (\xi_1, \xi_2, \ldots, \xi_n) \in \bR^d$ and for a.e. $(t,x) \in \Omega_T$.  We also assume that $a_0, c: \Omega_T \rightarrow \bR$ are given measurable functions satisfying
\begin{equation}  \label{a-b.zero0}
\nu \leq a_0(t,x), \ c(t,x) \leq \nu^{-1} \quad \text{for a.e.} \, (t,x) \in \Omega_T.
\end{equation}
We denote the following second-order linear operator in non-divergence form with singular coefficients
\begin{equation} \label{L.def}
\mathcal{L}u(t, x)= a_0(t, x) u_t  - a_{ij}(t, x) D_{ij}u  - \frac{\alpha}{x_d} a_{dj}(t, x) D_j u + \lambda c(t, x) u
\end{equation}
for $(t, x) = (t, x', x_d) \in \Omega_T$,
where $\alpha <1$ and $\lambda \geq 0$ are given. Our goal is to find a suitable class of Sobolev spaces for the well-posedness and regularity estimates of the following parabolic equations with the homogeneous Dirichlet boundary condition
\begin{equation} \label{main-eqn}
\left\{
\begin{array}{cccl}
\mathcal{L} u  & = & f & \quad \text{in} \quad \Omega_T,    \\
u & = & 0 & \quad \text{on} \quad (-\infty, T) \times \partial \bR^d_+.
\end{array} \right.
\end{equation}
When the coefficients $a_{ij}, c$, and $f$ are time independent, we also study the corresponding elliptic equations
\begin{equation} \label{elli-main-eqn}
\left\{
\begin{array}{cccl}
\sL u & = & f & \quad \text{in} \quad \bR^d_+,  \\
u & = & 0 & \quad \text{on} \quad \partial \bR^d_+,
\end{array}   \right.
\end{equation}
where
\begin{equation} \label{sL}
\sL u (x)= - a_{ij}( x) D_{ij}u   - \frac{\alpha}{x_d} a_{dj}(x) D_j u   +  \lambda c(x) u
\end{equation}
for $x =(x', x_d) \in \bR^d_+$. In Theorem \ref{para-main.theorem}, we show that under certain VMO conditions,  \eqref{main-eqn} has a unique solution in a suitable class of weighted mixed norm Sobolev spaces  with the weight $x_d^{p\alpha +\gamma_0}\omega_0(t)\omega_1(x)$ provided that $\lambda$ is sufficiently large. Here $\omega_0\in A_q(\bR)$ and $\omega_1\in A_p(\bR^d_+, \mu_2)$ are any Muckenhoupt weights for $q, p \in (1,\infty)$ and $\mu_2(s) = s^{\gamma_0}$ for $s>0$, and $\gamma_0 \in (-1, 1-\alpha]$. A similar result for the elliptic equation \eqref{elli-main-eqn} is stated in Theorem \ref{elli-main.theorem}.  From the mentioned theorems, we obtain the local boundary estimates stated in Corollary \ref{cor2.3}.  Under some mild conditions, it is possible that the mentioned results can be extended to the class of equations consisting a singular zeroth order term $bu/x_d^2$ as those considered in \cite{MNS21, MNS22}, by using a change of variables. See Remark \ref{final-remark} for details.

It should be mentioned that the estimates in our main results (Theorems \ref{thm0}, \ref{para-main.theorem}, and \ref{elli-main.theorem}) are quite different from those obtained in  \cite{Dong-Phan-RMI, DP} for the same class of equations but with the conormal boundary conditions, even when $p=2$. In fact, for a given solution $u$ of the PDE \eqref{main-eqn} or \eqref{elli-main-eqn},  $D_d^2 u$ could be too singular to be $L_p$-integrable  even with weights unless the weights have very fast decay near $\{x_d =0\}$. This can be seen by the ODE
\begin{equation} \label{ODE.exple}
u'' + \frac{\alpha}{x} u' =0 \quad \text{for} \,\, x \in (0,1)
\end{equation}
with a given $\alpha \in (0,1)$, for which $u(x) = x^{1-\alpha}$ is a solution and $u''(x) = -\alpha (1-\alpha) x^{-1-\alpha}$ which is strongly singular when $x\rightarrow 0^+$. This striking phenomenon can be seen clearly in the local pointwise estimates derived in Section \ref{simple-coeffients}. As such, instead of $D^2_d u$, we only derive the estimate for $D_d^2 u+\alpha D_du/x_d$. Therefore, in our main results, we establish the mixed-norm $L_p$-estimates of
$$
x_d^{\alpha} u,\,\,  x_d^{\alpha-1} u,
\,\, x_d^{\alpha} Du, \,\,  x_d^{\alpha-1} D_{x'}u,
\,\,
x_d^{\alpha} DD_{x'}u, \,\, x_d^{\alpha} u_t,\,\, \text{and}\,\,
x_d^{\alpha}(D_d^2 u+\alpha D_d u/x_d)
$$ with weight $\omega d\mu_2$ for a suitable nonnegative function $\omega$, while in \cite{Dong-Phan-RMI, DP} the mixed-norm $L_p$-estimates of $u,  Du,  D^2u, u_t$ with weight $\omega d\mu$ are obtained.  Due to such singularity feature for solutions of \eqref{main-eqn} and \eqref{elli-main-eqn}, suitable function spaces are found in Section \ref{func-setting}, and they are intrinsic for the problems \eqref{main-eqn} and \eqref{elli-main-eqn}.

The proofs of our main results are also different from those in \cite{Dong-Phan-RMI, Dong-Phan-1, Dong-Phan, DP}. More precisely, the proofs of the main results in  \cite{Dong-Phan-RMI, Dong-Phan-1, DP} use $\mu(s) = |s|^{\alpha}$ as an underlying measure, where $s \in \bR \setminus \{0\}$, and the proofs of the main results for equations in divergence form in \cite{Dong-Phan} use the underlying measure $\mu_1(s) =|s|^{-\alpha}$. In this paper, to prove Theorems \ref{para-main.theorem} - \ref{elli-main.theorem}, we introduce the new and more general underlying measure $\mu_2(s) = |s|^{\gamma_0}$ with $\gamma_0 \in (-1, 1-\alpha]$. Moreover, instead of the $L_2$-estimates as in \cite{Dong-Phan-1, DP}, the starting point in this paper is the weighted $L_p$-result stated in Lemma \ref{s-W-2-p.constant.eqn}, which is based on the weighted $L_p$ for divergence form equations established in \cite{Dong-Phan}.  By adapting the ideas in this paper, in certain cases it is possible to relax the condition \eqref{eq8.47} for non-divergence form equations with the boundary condition considered in \cite{DP}.

The remaining part of the paper is organized as follows. In the next section, we introduce the notation and the function spaces and state the main results of the paper. In Section \ref{Prelim}, we recall the definition of Muckenhoupt weights and state the weighted  mixed-norm Fefferman-Stein and Hardy-Littlewood maximal function theorems.
In Section \ref{simple-coeffients}, we consider equations with coefficients depending only on the $x_d$-variable. We first derive some local boundary estimates for  higher-order derivatives of solutions to homogeneous equations, which are the key estimates in the proof of the main theorems. In particular, we  prove Proposition \ref{Higher-Lip-est-lemma-s} below about pointwise estimates of solutions to homogeneous equations. Then, we use Proposition \ref{Higher-Lip-est-lemma-s} and an idea introduced in \cite{Kr99} to prove Theorem \ref{thrm-s6}, which is slightly more general than Theorem \ref{thm0}. Section \ref{VMO-section} is devoted to the proofs of Theorems \ref{para-main.theorem} and \ref{elli-main.theorem}, and Corollary \ref{cor2.3}. To prove Theorems \ref{para-main.theorem} and \ref{elli-main.theorem}, we apply the mean oscillation argument  in \cite{Krylov-book} with the underlying measure  $\mu_2$ mentioned above. Finally, to show Corollary \ref{cor2.3}, we use a localization and iteration argument.
%\newpage
\section{Function spaces, notation, and main results}     \label{sec2}

\subsection{Function spaces} \label{func-setting} For a given function $f$ defined in $\bR^d_+$ and for $\tau \in \bR$, we define the multiplicative operators
\[
\bM f (x) = x_d f(x) \quad \text{and} \quad \bM^\tau f(x) = x_d^\tau f(x) \quad \text{for } x = (x', x_d) \in \bR^d_+.
\]
Let $\sigma$ be a non-negative Borel measure on  either $\bR^d_+$ or $\bR^{d+1}_+$. For $p \in [1, \infty)$,  $-\infty\le S<T\le +\infty$, and $\cD \subset \bR^d_+$, and $Q:= (S,T)\times \cD$, let $L_p(Q, d\sigma)$ be the weighted Lebesgue space consisting of measurable functions $u$ on $Q$ such that the norm
\[
\|u\|_{L_p( Q, d\sigma)}= \bigg( \int_{Q} |u(t,x)|^p\, d\sigma (t,x) \bigg)^{1/p} <\infty.
\]
For $p,q\in [1,\infty)$, and the weights $\omega_0=\omega_0(t)$ and $\omega_1=\omega_1(x)$, we define $L_{q,p}(Q,\omega\,d\sigma)$ to be the weighted and mixed-norm Lebesgue space on $Q$ equipped with the norm
\begin{equation*}
				%			\label{eq0806_02}
\|u\|_{L_{q,p}((S,T)\times \cD, \omega\, d\sigma)}=\bigg(\int_{S}^T\Big(\int_{\cD} |u(t,x)|^p \omega_1(x)\,\sigma(dx)\Big)^{q/p}\omega_0(t)\,dt\bigg)^{1/q},
\end{equation*}
where $\omega(t,x)=\omega_0(t)\omega_1(x)$.  We define the weighted Sobolev space
$$
W^1_p(\cD,\omega_1\,d\sigma)=\big\{u\in L_p(\cD,\omega_1\,d\sigma):\,  Du\in L_p(\cD,\omega_1\,d\sigma)\big\}
$$
equipped with the norm
\[
\begin{split}
\|u\|_{W^1_p(\cD,\omega_1 d\sigma)}& =\|u\|_{L_p(\cD,\omega_1 d\sigma)} +\|Du\|_{L_p(\cD,\omega_1d\sigma)}.
\end{split}
\]
The Sobolev space $\sW^1_p(\cD, \omega_1 d\sigma)$ is defined to be the closure in $W^1_p(\cD,\omega_1\,d\sigma)$ of all compactly supported functions in $C^\infty(\overline{\cD})$ vanishing near $\overline{\cD} \cap \{x_d=0\}$.  We also define
\[
\sW^{2}_p(\cD,  \omega_1d\sigma) = \Big\{ u \in \sW^1_p(\cD, \omega_1 d\sigma): DD_{x'} u, \bM^{-\alpha}D_d(\bM^\alpha D_d u) \in L_{p}(\cD, \omega_1 d\sigma)\Big\},
\]
equipped with the norm
\[
\begin{split}
\|u\|_{\sW^{2}_p(\cD, \omega_1d\sigma)} = & \|u\|_{W^1_p(\cD,  \omega_1 d\sigma)} + \|DD_{x'}u\|_{ L_{p}(\cD, \omega_1 d\sigma)} \\
& \quad  + \|\bM^{-\alpha} D_d(\bM^\alpha D_d u)\|_{ L_{p}(\cD,  \omega_1 d\sigma)}.
\end{split}
\]
Similarly, for $Q = (S,T) \times \cD$, $\omega(t,x) = \omega_0(t) \omega_1(x)$, and for $q, p \in [1, \infty)$, we denote the mixed-norm weighted parabolic Sobolev space
\[
\begin{split}
\sW^{1,2}_{q,p}(Q, \omega d\sigma)  = &\Big \{u \in L_q((S,T), \sW^{2}_p(\cD,  \omega_1d\sigma), \omega_0),\  u_t \in L_{q,p}(Q, \omega d\sigma ) \Big\},
\end{split}
\]
equipped with the norm
\[
\begin{split}
&\|u\|_{\sW^{1,2}_{q,p}(Q,  \omega d\sigma)}  =  \left(\int_{S}^T \|u(t, \cdot)\|_{\sW^{2}_p(\cD, \omega_1d\sigma)}^q \omega_0(t) dt \right)^{1/q}  + \|u_t\|_{L_{q,p}(Q,  \omega d\sigma )}.
\end{split}
\]
We also denote $\hat{\sW}^{1,2}_{q,p}(Q,  \omega d\sigma)$  to be the subspace of $\sW^{1,2}_{q,p}(Q,  \omega d\sigma)$ defined by
\[
\hat{\sW}^{1,2}_{q,p}(Q,  \omega d\sigma) =\big\{ u \in \sW^{1,2}_{q,p}(Q,  \omega d\sigma): \bM^{-1}u, \bM^{-1} D_{x'}u \in L_{q,p}(Q,  \omega d\sigma ) \big\}
\]
and equipped with the norm
\[
\|u\|_{\hat{\sW}^{1,2}_{q,p}(Q,  \omega d\sigma)} = \|u\|_{\sW^{1,2}_{q,p}(Q,  \omega d\sigma)} + \|\bM^{-1} u\|_{L_{q,p}(Q,  \omega d\sigma )} + \|\bM^{-1} D_{x'}u\|_{L_{q,p}(Q,  \omega d\sigma )}.
\]
Similar, we also denote
\[
\hat{\sW}^2_p(\cD,  \omega_1d\sigma) = \big\{u \in \sW^{2}_p(\cD,  \omega_1d\sigma) : \bM^{-1} u, \bM^{-1} D_{x'}u \in L_p(\cD,  \omega_1d\sigma)\big\}
\]
equipped with the norm
\[
\|u\|_{\hat{\sW}^{2}_{p}(\cD,  \omega_1 d\sigma)} = \|u\|_{\sW^{2}_{p}(\cD,  \omega_1 d\sigma)} + \|\bM^{-1} u\|_{L_{p}(\cD,  \omega_1 d\sigma)} + \|\bM^{-1} D_{x'}u\|_{L_{p}(\cD,  \omega_1 d\sigma )}.
\]
The spaces $\hat{\sW}^{1,2}_{q,p}(Q,  x_d^{\alpha p}\omega d\sigma)$ and $\hat{\sW}^2_p(\cD,  x_d^{\alpha p}\omega_1d\sigma)$ are where the solutions of \eqref{main-eqn} and \eqref{elli-main-eqn} are found, respectively. However, in many intermediate steps, the results hold for  solutions in the larger spaces $\sW^{1,2}_{q,p}(Q,  \omega d\sigma)$ and $\sW^2_p(\cD, \omega_1d\sigma)$.

\begin{remark}
            \label{rem2.1}
It is clear that $\hat{\sW}^{1,2}_{q,p}(\Omega_T,\omega d\sigma)$ is a subspace of
\begin{align}
                            \label{eq3.18}
\tilde{\sW}^{1,2}_{q,p}(\Omega_T,\omega d\sigma)&:=\big\{u\in L_{q,p}(\Omega_T,\omega d\sigma):
u_t,\bM^{-1}u,Du,\bM^{-1}D_{x'}u,\notag\\
&\qquad DD_{x'}u,\bM^{-\alpha}D_d(\bM^\alpha D_d u) \in L_{q,p}(\Omega_T, \omega d\sigma)\big\}.
\end{align}
In fact, due to the term
$$
\|\bM^{-1} u\|_{L_{p}(\cD,  \omega_1 d\sigma)},
$$
it follows from Lemma \ref{lemma-3.4} below that these two spaces are identical.
\end{remark}

\subsection{Notation and main results} Let $r>0$, $z_0 = (t_0, x_0)$ with $x_0 = (x_0', x_{0d}) \in  \bR^{d-1} \times \bR$ and $t_0 \in \bR$. We define $B_r(x_0)$  to be the ball in $\bR^d$ of radius $r$ centered at $x_0$, $Q_r(z_0)$ to be the parabolic cylinder of radius $r$ centered at $z_0$:
\[
Q_r(z_0) = (t_0-r^2 , t_0)\times B_r(x_0).
\]
Also, let $B_r^+(x_0)$ and $Q_r^+(z_0)$ be the upper-half ball and cylinder of radius $r$ centered at $x_0$ and $z_0$, respectively:
\begin{align*}
B_r^+(x_0) &= \big\{x = (x', x_d) \in \bR^{d-1} \times \bR:\, x_d >0, \ |x -x_0| < r\big\},\\
Q_r^+ (z_0)&= ( t_0-r^2, t_0)\times B_r^+ (x_0).
\end{align*}
For $z_0' = (t_0, x_0') \in \bR \times \bR^{d-1}$, we denote the parabolic cylinder in $\bR \times \bR^{d-1}$ by
$$Q'_\rho(z_0') = (t_0 -\rho^2, t_0) \times B_\rho'(x_0'),$$
where $B_\rho'(x_0')$ is the ball in $\bR^{d-1}$ of radius $\rho$ centered at $x_0'$.   Throughout the paper, when $x_0 =0$ and $t_0 =0$, for  simplicity of notation, we drop $x_0, z_0$ and write $B_r$, $B_r^+$, $Q_r$, and $Q_r^+$, etc.

For a measurable set $\Omega\subset \bR^{d+1}$ and any integrable function $f$ on $\Omega$ with respect to some locally finite Borel measure $\sigma$, we write
$$
\fint_\Omega f(z)\ \sigma(dz) =\frac 1 {\sigma(\Omega)}\int_\Omega f(z)\, \sigma(dz), \quad \text{where } \sigma(\Omega) = \int_{\Omega} \sigma(dz).
$$
Throughout the paper, for $\alpha \in (-\infty, 1)$ and a number $\gamma_0 \in (-1, 1-\alpha]$, we denote the following weights which are used frequently in the paper
\begin{equation} \label{mu-def}
 \mu(s) = |s|^\alpha, \quad \mu_1(s) = |s|^{-\alpha}, \quad \mu_2(s) =|s|^{\gamma_0} \quad \text{for } s \in \bR \setminus\{0\}.
\end{equation}
We  also write
\[
\begin{split}
& \mu(dz) = \mu(x_d) \, dx dt, \quad  \mu(dx) = \mu(x_d)\, dx, \\
& \mu_k(dz) = \mu_k(x_d) \, dx dt,  \quad\mu_k(dx) = \mu_k(x_d)\, dx, \quad k =1, 2.
 \end{split}
\]
For any $z_0 = (z_0' , x_{d0}) \in \overline{\Omega}_T$ and $\rho >0$,  we define the average of $f$ in $Q'_\rho(z_0')$ as
\begin{equation} \label{average-a}
 [f]_{\rho, z_0}(x_d)  =\fint_{Q'_\rho(z_0')} f(t, x', x_d)\, dx'dt
\end{equation}
and its weighted average in $Q_\rho^+(z_0)$ as
\begin{equation} \label{k-mean}
(f)_{\rho,z_0}=\fint_{Q^+_\rho(z_0)} f(z) \ \mu_2(dz).
\end{equation}
We denote the weighted mean oscillation of the given coefficients $(a_{ij})$, $a_0$, and $c$ by  \begin{align} \label{a-sharp}
a_{\rho}^\#(z_0) &= \sum_{i=1}^{d-1} \sum_{j=1}^d\fint_{Q_\rho^+(z_0)}\Big|\tilde{a}_{ij}(z)-[\tilde{a}_{ij}]_{\rho, z_0}(x_d)\Big|\,\mu_2(dz)\nonumber\\
& \quad + \sum_{j=1}^{d-1}\fint_{Q_\rho^+(z_0)}\Big|\tilde{a}_{dj}(z)-(\tilde{a}_{dj})_{\rho, z_0}\Big|\,\mu_2(dz)\nonumber\\
&\quad +\fint_{Q_\rho^+(z_0)}\Big(\big|\tilde{a}_0(z)-[\tilde{a}_0]_{\rho, z_0}(x_d)\big|+\big|\tilde{c}(z)-[\tilde{c}]_{\rho, z_0}(x_d)\big|\Big)\,\mu_2(dz)
\end{align}
for $z_0\in \overline{\Omega_T}$, where
\[
 \tilde{a}_0 = a_0/a_{dd}, \quad \tilde{c} = c/a_{dd}, \quad \text{and} \quad \tilde{a}_{ij} = a_{ij}/a_{dd}, \quad \text{for} \quad i, j =1,2,\ldots, d.
\]
When the coefficients are time-independent, we similarly define $a_{\rho}^\#(x_0)$ for $x_0 \in  \overline{\bR_+^d}$. We point out that in \cite{Dong-Phan} where the corresponding class of \eqref{main-eqn} in divergence form is considered, the mean oscillations of the coefficients are measured with the weight $\mu_1$.  Here we use $\mu_2$ in \eqref{a-sharp} and this is more general as when $\gamma_0 = -\alpha$, we have $\mu_2 \equiv \mu_1$. See the work \cite{Dong-Phan-RMI, Dong-Phan-1, DP, MP} for similar definitions of mean oscillations but with different weights.

\begin{remark}
We note that in the definition of $a_{\rho}^\#(z_0)$, the mean oscillations of $\tilde{a}_{dj}$ with $j=1,2,\ldots, d-1$ are measured in all variables. For the other coefficients, their mean oscillations are measured only in $(t,x')$. A smallness condition on such partial mean oscillations of the coefficients was introduced in \cite{Kim-Krylov,Kim-Krylov-1}. It is clear that with the weight, the smallness condition is weaker for larger $\gamma_0$. Also, as $\tilde{a}_{dd} \equiv 1$, its mean oscillation is zero hence it does not appear in $a_{\rho}^\#(z_0)$.  Observe also that we may assume without loss of generality that $a_{dd} \equiv 1$ as we can always divide both sides of the PDE in \eqref{main-eqn} by $a_{dd}$ and replace $\nu$ in \eqref{ellipticity} and \eqref{a-b.zero0} with $\nu^2$.
\end{remark}

By a strong solution $u \in \hat{\sW}^{1,2}_{q,p}(\Omega_T,  x_d^{\alpha p} \omega\, d\sigma)$ to \eqref{main-eqn} with $p, q \in (1, \infty)$, we mean that the first equation of \eqref{main-eqn} is satisfied almost everywhere. By a strong solution $u\in \sW^{1,2}_{q,p}(\Omega_T,  x_d^{\gamma} dz)$ to \eqref{main-eqn}, we mean that the first equation of \eqref{main-eqn} is satisfied almost everywhere and the zero Dirichlet boundary condition is satisfied in the sense of trace.
Note that the solution space $\hat{\sW}^{1,2}_{q,p}(\Omega_T,  x_d^{\alpha p}\omega\, d\sigma)$ (or $\sW^{1,2}_{q,p}(\Omega_T,  x_d^{\gamma} dz)$) is included in the usual parabolic Sobolev space $W^{1,2}_{q,p,\text{loc}}(\Omega_T, \omega d\sigma)$  (or $\sW^{1,2}_{q,p,\text{loc}}(\Omega_T,  dz)$, respectively), so that the derivatives of $u$ on the left-hand side of \eqref{main-eqn} are defined almost everywhere.  Moreover, the trace operator is well defined for $\sW^{1,2}_{q,p}(\Omega_T,  x_d^{\gamma}dz)$ when $\gamma<2p-1$.

We are now ready to state the first main result of the paper.
\begin{theorem}
\label{para-main.theorem}  Let $\nu \in (0,1)$, $T\in (-\infty,\infty]$, $p, q, K \in (1, \infty)$, $\alpha \in (-\infty,1), \gamma_0 \in (-1, 1-\alpha]$, and $\rho_0 >0$. Then
there exist $\delta = \delta (d, \nu, p, q, \alpha, \gamma_0, K) >0$ sufficiently small and $ \lambda_0 =  \lambda_0 (d,\nu,p,q,\alpha, \gamma_0, K)>0$ such that the following assertion holds. Suppose that \eqref{ellipticity} and \eqref{a-b.zero0} are satisfied, $\omega_0\in A_q(\bR)$, $\omega_1 \in A_p(\bR^{d}_+, \mu_{2})$  with
$$
[\omega_0]_{A_q(\bR)}, \quad  [\omega_1]_{A_p(\bR^d_+, \mu_{2})} \le K,
$$
and
\begin{equation} \label{para-VMO}
a^{\#}_{\rho}(z_0) \leq \delta, \quad \forall \ \rho \in (0, \rho_0), \quad \forall \  z_0 \in  \overline{\Omega}_T.
\end{equation}
Then for any $f \in L_{q,p}(\Omega_T, x_d^{p\alpha}\omega\, d\mu_2)$ and $\lambda \ge  \lambda_0 \rho_0^{-2}$, there exists a unique strong solution $u \in \hat{\sW}^{1,2}_{q,p}(\Omega_T,  x_d^{p\alpha}{\omega}\, d\mu_2)$ to \eqref{main-eqn}, which satisfies
\begin{equation} \label{main-para.est}
\begin{split}
& \norm{u_t}_{L_{q,p}} + \norm{DD_{x'}u}_{L_{q,p}}
+ \norm{\bM^{-\alpha} D_d(\bM^\alpha D_d u)}_{L_{q,p}} +\norm{\bM^{-1} D_{x'}u}_{L_{q,p}} \\
&\quad + \sqrt{\lambda} \norm{Du}_{L_{q,p}} + \lambda \norm{u}_{L_{q, p}} +\sqrt\lambda \norm{\bM^{-1} u}_{L_{q, p}} \leq  N\norm{f}_{L_{q,p}},
\end{split}
\end{equation}
where $\omega(t,x) = \omega_0(t) \omega_1(x)$ for $(t,x) \in \Omega_T$,   $d\mu_2 = x_d^{\gamma_0} \,dxdt$,
\[ L_{q, p} = L_{q,p}(\Omega_T, x_d^{p\alpha}\omega\, d\mu_2), \quad \text{and} \quad N = N( d,\nu, p, q, \alpha, \gamma_0, K) >0. \]
\end{theorem}

For elliptic equations, we also obtain the following results concerning \eqref{elli-main-eqn}.
\begin{theorem} \label{elli-main.theorem} Let $\nu \in (0,1), p,K \in (1, \infty)$,  $\alpha \in (-\infty,1), \gamma_0 \in (-1, 1-\alpha]$, and $\rho_0>0$. There exist $\delta = \delta (d, \nu, p, \alpha, \gamma_0, K) >0$ sufficiently small and $\lambda_0 = \lambda_0 (d,\nu,p,\alpha, \gamma_0, K)>0$ such that the following assertion holds.  Suppose that \eqref{ellipticity}, and \eqref{a-b.zero0} are satisfied, $\omega \in A_p(\bR_+^d, \mu_{2})$ with $[\omega]_{A_p(\bR_+^d, \mu_2)} \le K$, and
\begin{equation} \label{elli-VMO}
a^{\#}_{\rho}(x_0) \leq \delta, \quad \forall \ \rho \in (0, \rho_0), \quad \forall  \ x_0 \in  \overline{\bR_+^d}.
\end{equation}
Then for any $f \in L_{p}(\bR^d_+, x_d^{p\alpha}\omega\, d\mu_2)$ and for $\lambda \geq \lambda_0\rho_0^{-2}$, there exists a unique strong solution $u \in \hat{\sW}^{2}_{p}(\bR^d_+, x_d^{p\alpha}\omega\, d\mu_2)$ to \eqref{elli-main-eqn}, which satisfies
\begin{equation} \label{main-ell.est}
\begin{split}
&  \norm{DD_{x'}u}_{L_{p}}  +  \norm{\bM^{-\alpha} D_d(\bM^\alpha D_{d}u)}_{L_{p}} +\norm{\bM^{-1} D_{x'}u}_{L_{p}} \\
&\quad  + \sqrt{\lambda} \norm{Du}_{L_{p}} + \lambda \norm{u}_{L_{p}}  +\sqrt\lambda \norm{\bM^{-1} u}_{L_{ p}} \leq  N \norm{f}_{L_{p}},
 \end{split}
\end{equation}
where  $L_p=L_{p}(\bR_+^d, x_d^{p\alpha}\omega\, d\mu_2)$, $N = N( d,\nu, p, \alpha, \gamma_0, K) >0$ and $d\mu_2 = x_d^{\gamma_0} \, dx$.
\end{theorem}

 A few remarks about the theorems above are in order.
\begin{remark} \label{rem2.4}
A typical example of weights is the power weights $\omega_{1}(x_d)=x_d^\beta$. It is easily seen that $\omega_{1} \in A_p(\bR^{d}_+, \mu_2)$ if and only if $\beta \in (-\gamma_0 -1, (1+\gamma_0) (p-1))$. Therefore, from Theorem \ref{para-main.theorem}, we obtained the estimate and solvability in the space $\hat{\sW}^{1,2}_{q,p}(\Omega_T,x_d^{\gamma} dz)$, where $\gamma=\beta+\alpha p+\gamma_0\subset (\alpha p -1, (1+\alpha +\gamma_0)p-1)$. In the special case when $\alpha=0$, similar results were obtained in \cite{K99,KKL,DK15}. However, the powers of the distance function in these papers vary with
the order of derivatives. Thus the results in these papers cannot be directly deduced from
Theorem \ref{para-main.theorem}.
\end{remark}
%===
\begin{remark} \label{remark2}
Theorems \ref{para-main.theorem}-\ref{elli-main.theorem}  imply Theorem \ref{thm0} in the introduction. In fact, when the coefficients $a_{ij}, a_0, c$ are constants, the conditions \eqref{para-VMO} and \eqref{elli-VMO} hold for all $\gamma_0$ and for all $\rho_0>0$. Then, as $\gamma \in (\alpha p -1, 2p -1)$, we can choose $\gamma_0 \in (-1, 1-\alpha]$ such that $\gamma=\beta+\alpha p + \gamma_0$ with $\beta \in (-\gamma_0 -1, (1+\gamma_0) (p-1))$. From this, Remark \ref{rem2.4}, and a standard scaling argument $u(t,x)\to u(s^2 t,sx)$ for $s>0$, we see that \eqref{main-para.est} and \eqref{main-ell.est} hold for any $\lambda > 0$. See Theorem \ref{thrm-s6} below for a slightly generalization of Theorem \ref{thm0}. When $\gamma \in ((\alpha+1)p-1, 2p -1)$, we also obtain the estimates of $\|D_d^2 u\|_{L_p(\Omega_T, x_d^\gamma dz)}$, $\|\bM^{-1} Du\|_{L_p(\Omega_T, x_d^{\gamma} dz)}$, and $\|\bM^{-2} u\|_{L_p(\Omega_T, x_d^{\gamma} dz)}$, When $\alpha=0$, this agrees with the results in \cite{K99,DK15}. See Lemma \ref{s-W-2-p.constant.eqn} and Remark \ref{rm-47} below.
\end{remark}

Finally, we state a local estimate, which is a consequence of Theorems \ref{para-main.theorem}--\ref{elli-main.theorem}.

\begin{corollary}
                    \label{cor2.3}
Let $\nu \in (0,1)$, $p, q, K \in (1, \infty)$, $\alpha \in (-\infty,1), \gamma_0 \in (-1, 1-\alpha]$, $\lambda\in [0,\infty)$, and $\rho_0 >0$. Then
there exists $\delta = \delta (d, \nu, p, q, \alpha, \gamma_0, K) >0$ sufficiently small such that the following assertion holds. Suppose that \eqref{ellipticity}, \eqref{a-b.zero0}, and \eqref{para-VMO} are satisfied, $\omega_0\in A_q(\bR)$, $\omega_1 \in  A_p(\bR^{d}_+, \mu_2)$  with
$$
[\omega_0]_{A_q(\bR)}, \,\,  [\omega_1]_{A_p(\bR^{d}_+, \mu_2)}  \le K.
$$
Assume that $f \in L_{q,p}(Q_1^+, x_d^{p\alpha}\omega\, d\mu_2)$ and $u \in \hat{\sW}^{1,2}_{q,p}(Q_1,  x_d^{p\alpha}{\omega}\, d\mu_2)$ is strong solution of \eqref{main-eqn} in $Q_1^+$. Then we have
\begin{equation} \label{eq2.17}
\begin{split}
& \norm{u_t}_{L_{q,p}(Q_{1/2}^+,x_d^{p\alpha}\omega\, d\mu_2)} + \norm{DD_{x'}u}_{L_{q,p}(Q_{1/2}^+,x_d^{p\alpha}\omega\, d\mu_2)}\\
&\quad + \norm{D_d(\bM^\alpha D_d u)}_{L_{q,p}(Q_{1/2}^+, \omega d\mu_2)}  + \norm{\bM^{-1} D_{x'}u}_{L_{q,p}(Q_{1/2}^+,x_d^{p\alpha}\omega\, d\mu_2)}\\
 &\quad +\norm{Du}_{L_{q,p}(Q_{1/2}^+,x_d^{p\alpha}\omega\, d\mu_2)}
+\norm{\bM^{-1} u}_{L_{q,p}(Q_{1/2}^+,x_d^{p\alpha}\omega\, d\mu_2)}
\\
& \leq  N\norm{f}_{L_{q,p}(Q_{1}^+,x_d^{p\alpha}\omega\, d\mu_2)}+\norm{u}_{L_{q, p}(Q_{1}^+,x_d^{p\alpha}\omega\, d\mu_2)},
\end{split}
\end{equation}
where $\omega(t,x) = \omega_0(t) \omega_2 (x)$ for $(t,x)\in Q_1^+$, $N = N(\nu, d, p, q, \alpha, \gamma_0, K) >0$, and $d\mu_2 = x_d^{\gamma_0} \, dxdt$. A similar local estimate holds for the elliptic equation \eqref{elli-main-eqn} in $B_1^+$.
\end{corollary}
\begin{remark} It is worth pointing out that \eqref{eq8.47} is required in \cite{Dong-Phan-RMI, DP}. Similar structural conditions on the matrix $(a_{ij})$ are also imposed in \cite{MP, Sire-1,Sire-2}. In this paper, we do not have those restrictions. This is due to  the new H\"{o}lder{\color{red}'s} regularity of $x_d^{\alpha-1}D_{x'}^ku$ established in \eqref{eq4.39} in Lemma \ref{Higher-Lip-est-lemma} below when $u$ solves the homogeneous equations.
\end{remark}

\section{Preliminaries on weights and weighted inequalities} \label{Prelim}

We first recall the definition of Muckenhoupt weights, which was introduced in \cite{Muck}.
\begin{definition}
        \label{Def-Muck-wei}
For each $p \in (1, \infty)$  and for a non-negative Borel measure $\sigma$ on $\bR^d$, a locally integrable function $\omega :  \bR^d \rightarrow \bR_+$ is said to be in the $A_p( \bR^d,  \sigma)$ Muckenhoupt class of weights if and only if $[\omega]_{A_p(\bR^d, \sigma)} < \infty$, where
\begin{equation}
                    \label{Ap.def}
[\omega]_{A_p(\bR^d,  \sigma)} =
\sup_{\rho >0,x \in \bR^d } \bigg[\fint_{B_\rho (x)} \omega(y)\,  \sigma(dy) \bigg]\bigg[\fint_{B_\rho(x)} \omega(y)^{\frac{1}{1-p}}\,  \sigma(dy) \bigg]^{p-1}.
\end{equation}
Similarly, the class of weight $A_p(\bR^d_+,  \sigma)$ can be defined in the same way in which the ball $B_\rho (x)$ in \eqref{Ap.def} is replaced with $B_\rho^+(x)$ for $x\in \overline{\bR^d_+}$. If $\sigma$ is a Lebesgue measure, we simply write $A_p(\bR^d_+) = A_p(\bR^d_+, dx)$ and $A_p(\bR^d) = A_p(\bR^d, dx)$.  Note that if $\omega \in A_p(\bR)$, then $\tilde{\omega} \in A_p(\bR^d)$ with $[\omega]_{A_p(\bR)} = [\tilde{\omega}]_{A_p(\bR^d)}$, where $\tilde{\omega}(x) = \omega(x_n)$ for $x = (x', x_n) \in \bR^d$. Sometimes, if the context is clear, we neglect the spatial domain and only write $\omega \in A_p$.
\end{definition}

Denote the collection of parabolic cylinders in $\Omega_T$ by
\[
\mathcal{Q} = \{ Q_\rho^+(z): \rho >0, z \in \Omega_T \}.
\]
Now, for any locally integrable function $f$ defined in $\Omega_T$, the Hardy-Littlewood maximal function of  $f$ with respect to $d\mu_2$ is defined by
\begin{equation} \label{maxi-def}
\mathcal{M}(f)(z) = \sup_{Q \in \mathcal{Q}, z \in Q}\fint_{Q} |f(\xi)| \ \mu_2(d\xi),
\end{equation}
and the Fefferman-Stein sharp function of $f$ with respect to $d\mu_2$ is defined by
\begin{equation} \label{sharp-def}
f^{\#}(z) =  \sup_{Q \in \mathcal{Q},z  \in Q} \fint_{Q}|f(\xi) - (f)_{Q}| \ \mu_2(d\xi),
\end{equation}
where $\mu_2$ is defined in \eqref{mu-def}, and
\begin{equation} \label{k-meanb}
(f)_{Q}=\fint_{Q} f(z) \ \mu_2(dz).
\end{equation}

The following version of the weighted mixed-norm Fefferman-Stein theorem and Hardy-Littlewood maximal function theorem can be found in \cite{Dong-Kim-18}.
\begin{theorem}  \label{FS-thm} Let $p, q \in (1,\infty)$, $\gamma_0 \in (-1, 1-\alpha]$, $K\geq 1$. Suppose that $\omega_0\in A_q(\bR)$, $\omega_1 \in A_p(\bR^{d}_{+},\mu_2)$ with
$$
[\omega_0]_{A_q},    \,\, [\omega_{1}]_{A_p(\bR_+^d, \mu_2)}\le K.$$
Then, for any $f \in L_{q, p}(\Omega_T, \omega\, d\mu_2)$, we have
\begin{equation*} %\label{Maximal-L-p}
\begin{split}
& \|f\|_{L_{q, p}(\Omega_T,  \omega\, d\mu_2)} \leq N \| f^{\#}\|_{L_{q,p}(\Omega_T,  \omega\, d \mu_2)} \quad \text{and} \quad \\
& \|\mathcal{M}(f)\|_{L_{q,p}(\Omega_T, \omega\, d\mu_2)} \leq N \|f\|_{L_{q, p}(\Omega_T,  \omega \, d \mu_2)},
\end{split}
\end{equation*}
where $N = N(d, q, p, \gamma_0, K)>0$ and $\omega(t,x) = \omega_0(t)\omega_1(x)$ for $(t,x) \in \Omega_T$.
\end{theorem}
We now state the well-known weighted Hardy's inequalities whose proof can be found, for instance, in \cite[Lemma 3.1]{Dong-Phan}.
\begin{lemma}\label{w-hardy-inq} For $p \in [1, \infty)$, the following statements hold.
\begin{itemize}
\item[\textup{(i)}] For each $\beta +1 < p$ and a measurable function $f$ defined on $\bR_+$, we have
\[
\|g\|_{L_p(\bR_+, s^{\beta-p}ds)} \leq \frac{p}{p-(1+\beta)} \|f\|_{L_p(\bR_+, s^{\beta} ds)},
\]
where $g(s) = \int_0^s f(\tau)\, d\tau$.
\item[\textup{(ii)}] For each $\beta +1 > p$ and measurable function $f$ defined on $\bR_+$, we have
\[
\|g\|_{L_p(\bR_+, s^{\beta-p}ds)} \leq \frac{p}{\beta+1-p} \|f\|_{L_p(\bR_+, s^{ \beta} ds)},
\]
where $g(s) = \int_s^\infty f(\tau)\, d\tau$.
\end{itemize}
\end{lemma}
%===

\begin{lemma}
        \label{lemma-3.4}
Let $p \in [1, \infty)$,  $\omega$ be a weight, and $\sigma$ be a locally finite non-negative Borel measure on $\bR^d_{+}$  such that the set of continuous functions in $\bR^d_+$ is dense in $L_p(\bR^d_+,\omega d\sigma)$. Assume that  $u \in W^{1}_p(\bR^d_+, \omega d\sigma)$ and
\[
\|\bM^{-1} u\|_{L_p(\bR^d_+, \omega d\sigma)} < \infty.
\]
Then, there exists a sequence of smooth functions $\{u_k\}$ in $W^1_p(\bR^d_+, \omega d\sigma)$ vanishing near $\{x_d =0\}$,  which converges to $u$ in $W^1_p(\bR^d_+, \omega d\sigma)$.
\end{lemma}
\begin{proof} Let $\eta \in C^\infty(\bR)$ be such that $\eta(s) =0$ for $s \leq  1/2$ and $\eta(s) =1$ for $s \geq 1$. For each $k \in \mathbb{N}$, let $v_k(x) = u(x) \eta_k(x_d)$ for $x = (x', x_d) \in \bR^d_+$, where
\[ \eta_k (s) = \eta(ks), \quad s \in \bR.\]
By the Lebesgue dominated convergence theorem, we see that for $j=1,2,\ldots, d-1$,
\[
v_k \rightarrow u \quad \text{and} \quad  \quad D_{x'} v_k \rightarrow D_{x'} u
\]
in $L_p(\bR^d_+, \omega d\sigma)$ as $k\rightarrow \infty$. Now note that
\[
D_d v_k(x) = \eta_k(x_d) D_d u(x) + k \eta'(kx_d)  u(x).
\]
Since
\[
|k \eta'(kx_d)| \leq N x_d^{-1}{\bf 1}_{(0, 1/k)}(x_d),
\]
as $k\rightarrow \infty$,
\[
\left(\int_{\bR^{d}_+} |k \eta'(kx_d)  u(x)|^p \omega(x) d\sigma(x)\right)^{1/p} \leq N \|\bM^{-1} u\|_{L_p(\bR^{d-1} \times (0, 1/k), \omega\, d\sigma)} \rightarrow 0.
\]
From this and by using the dominated convergence theorem, we obtain
\[
D_d v_k \rightarrow D_d u \quad \text{in} \quad L_p(\bR^d_+, \omega\, d\sigma) \quad \text{as} \quad k \rightarrow \infty.
\]
Consequently, $\{v_k\}$ converges to $u$ in $W^1_p(\bR^d_+, \omega d\sigma)$ and $v_k$ vanishes near $\{x_d =0\}$ for each $k$. Finally, by using  the standard mollification, we can  find a sequence of smooth functions $\{u_k\}$ in $W^1_p(\bR^d_+, \omega d\sigma)$ satisfying the assertion of the lemma.
\end{proof}
%===
%\todo{Regarding your concern, to have a trace operator, we need $p>\gamma+1$, so $x_d^\beta\notin W^1_p(x_d^\gamma)$ when $\beta<0$.}
We conclude the section with the following lemma, which is used frequently in the paper.
\begin{lemma}
            \label{D-dd-control}
Let $\nu \in (0,1), \alpha \in (-\infty,1)$ and $p ,q \in (1, \infty)$.  Let $\sigma$ be a non-negative Borel measure $\bR^{d+1}_+$ and $\omega : \Omega_T \rightarrow \bR_+$ be a weight. Suppose that \eqref{ellipticity} and \eqref{a-b.zero0} are satisfied.  Then for
any $R \in (0, \infty]$, if $u$ is a strong solution of
\[
\mathcal{L} u   =  f   \quad \text{in} \quad Q_{R}^+
 \]
with some $\lambda \ge 0$ and $f \in L_{q, p}(Q_R^+, x_d^{\alpha p} \omega\, d\sigma)$, then it holds that
\begin{equation*} % \label{Ddd-u.est}
\begin{split}
& \|D_{d}(\bM^\alpha D_d u)\|_{L_{q,p}(Q_R^+,  \omega\, d\sigma)}  \leq  N \Big[\|u_t\|_{L_{q,p}(Q_R^+,  x_d^{\alpha p}\omega\, d\sigma)} + \|DD_{x'}u\|_{L_{q,p}(Q_R^+, x_d^{\alpha p} \omega\, d\sigma)} \\
 & \qquad + \|\bM^{-1}D_{x'}u\|_{L_{q,p}(Q_R^+, x_d^{\alpha p} \omega\, d\sigma)}  +  \lambda \|u\|_{L_{q,p}(Q_R^+, x_d^{\alpha p} \omega\, d\sigma)} + \|f\|_{L_{q,p}(Q_R^+, x_d^{\alpha p} \omega\, d\sigma)} \Big],
\end{split}
\end{equation*}
where $N= N(d,\nu, \alpha) >0$.
\end{lemma}
\begin{proof} Without loss of generality, we may assume that the right-hand side of the inequality above is finite.  By dividing the PDE of $u$ by $a_{dd}$ and using  the conditions \eqref{ellipticity} and \eqref{a-b.zero0}, we obtain
\[
|D_{d}(\bM^{\alpha} D_du)| \leq N(d, \nu, \alpha) \bM^{\alpha} F,
\]
where
\[
F= |f| + \lambda |u| + \bM^{-1}|D_{x'}u| + |u_t| + |DD_{x'}u|.
\]
Therefore,
\[
 \|D_{d}(\bM^\alpha D_d u)\|_{L_{q,p}(Q_R^+,  \omega\, d\sigma)}\leq  N \|F\|_{L_{q,p}(Q_R^+,  x_d^{\alpha p}\omega\, d\sigma)}.
\]
The lemma is proved.
\end{proof}
%=====
\section{Equations with simple coefficients} \label{simple-coeffients}

We consider the special class of equations \eqref{main-eqn} in which the coefficients only depend on the $x_d$-variable. Let $(\overline{a}_{ij}) : \bR_+ \rightarrow \bR^{d\times d}$ be bounded, measurable, and uniformly elliptic: there is $\nu \in (0,1)$ so that
\begin{equation} \label{con-elli}
\nu |\xi|^2 \leq \overline{a}_{ij}(x_d) \xi_i \xi_j \quad \text{and} \quad |\overline{a}_{ij}(x_d)| \leq \nu^{-1}
\end{equation}
for $x_d \in \bR_+$ and for $\xi = (\xi_1,\xi_2,\ldots, \xi_d) \in \bR^d$. Moreover, let $\overline{a}_0, \overline{c}: \bR_+ \rightarrow \bR$ be measurable functions satisfying
\begin{equation}  \label{a-b.zero}
\nu \leq \overline{a}_0(x_d), \ \overline{c}(x_d) \leq \nu^{-1} \quad \text{for a.e. } x_d \in \bR_+.
\end{equation}
For each $\alpha <1$ and $\lambda \geq 0$, we denote
\begin{equation*} %\label{L-0.def}
\mathcal{L}_0 u(t, x) = \overline{a}_0(x_d) u_t + \lambda  \overline{c}(x_d) u -  \overline{a}_{ij}(x_d) D_{ij} u(t, x', x_d) - \frac{\alpha}{x_d} \overline{a}_{dj} D_j u(t, x', x_d),
\end{equation*}
where $(t, x)  = (t, x', x_d) \in \Omega_T$. We consider the following equation
\begin{equation} \label{constant-L-0}
\left\{
\begin{array}{ccll}
  \mathcal{L}_0 u  & = & f & \quad  \text{in} \quad \Omega_T, \\
u(\cdot,0) & = & 0 & \quad  \text{on} \quad  (-\infty, T) \times \bR^{d-1}.
\end{array}  \right.
\end{equation}
In addition to the uniformly elliptic and bounded conditions as in \eqref{con-elli}, we assume that
\begin{equation}
                            \label{eq4.31}
\overline{a}_{dj}/\overline{a}_{dd}, \quad j = 1, 2,\ldots, d-1\text{ are constant.}
\end{equation}
Dividing both sides of the equation by $\overline{a}_{dd}$, we may assume that
\begin{equation} \label{a-dd.condb}
\overline{a}_{dj}(x_d) \equiv \overline{a}_{dj} \quad \text{and} \quad \overline{a}_{dd}(x_d) \equiv 1, \quad \forall \ x_d \in \bR_+, \quad j = 1, 2,\ldots, d-1.
\end{equation}
Observe that under this assumption and by a change of variables, $y_j=x_j- \overline{a}_{dj}x_d,j=1,2,\ldots,d-1$ and $y_d=x_d$, without loss of generality, we may assume that $\overline{a}_{dj} \equiv 0$ for $j = 1,2,\ldots, d-1$ as in \eqref{extension-type-matrix}.  Hence, in the remaining part of this section, we assume that
\begin{equation} \label{a-dd.cond}
\overline{a}_{dj}(x_d) \equiv  0 \quad \text{and} \quad \overline{a}_{dd}(x_d) \equiv 1, \quad \forall \ x_d \in \bR_+, \quad j = 1, 2,\ldots, d-1.
\end{equation}
We stress that the results below still hold true under the condition \eqref{eq4.31} by changing the variables back.

Observe that under the condition \eqref{a-dd.cond}  or \eqref{a-dd.condb}, there is a hidden divergence structure for the operator $\mathcal{L}_0$. Namely,
\[
x_d^\alpha \mathcal{L}_0 u (t, x) =x_d^\alpha\big( \overline{a}_0(x_d) u_t + \lambda  \overline{c}(x_d) u\big) - D_i[x_d^\alpha \overline{a}_{ij}(x_d) D_j u(t, x) ].
\]
Consequently, the PDE in \eqref{constant-L-0} can be rewritten
in divergence form as
\begin{equation} \label{div-L-0}
x_d^{\alpha}\big(\overline{a}_0(x_d)u_t + \lambda \overline{c}(x_d) u \big)- D_i[x_d^\alpha \overline{a}_{ij}(x_d) D_j u(t,x)]   =  x_d^\alpha f(t,x)  \quad \text{in} \quad \Omega_T.
\end{equation}
A function $u \in L^2((-\infty,T), \sW^1_p(\bR^d_+, d\mu))$ is said to be a weak solution of \eqref{constant-L-0} if
\begin{align*}  %\label{weak-formula-bdr}
\int_{\Omega_T}\mu(x)[- \overline{a}_0 u\varphi_t + \overline{a}_{ij}D_j u D_i\varphi+\lambda \overline{c} u\varphi ]\,dz = \int_{\Omega_T} \mu(x)f \varphi \,dz
\end{align*}
for any $\varphi \in C_0^\infty(\Omega_T)$ and for $\mu(x) = x_d^\alpha$ with $x = (x', x_d) \in \bR^{d}_+$.
\subsection{Local pointwise estimates for homogeneous equations}  We consider the equation
\begin{equation} \label{int-Q4-constant-L-0}
\left \{
\begin{array}{ccll}
 \mathcal{L}_0 u  &  = &  0 &\quad  \text{in}\quad Q_{2}^+(\hat{z}) \\
 u & = & 0 & \quad \text{on} \quad Q_{2}(\hat{z})\cap \{x_d =0\} \quad \text{if} \quad \hat{x}_d \leq 2,
\end{array} \right.
\end{equation}
where $\hat{z} = (\hat{t}, \hat{x}', \hat{x}_d) \in \bR \times \bR^{d-1} \times \overline{\bR_+}$.
Our goal is to derive pointwise estimates for solutions to \eqref{int-Q4-constant-L-0} and their derivatives. We start with the following Caccioppoli type estimates.
\begin{lemma}
            \label{cappio}
Let $\nu \in (0, 1]$, $\lambda \geq 0$,  $\alpha <1$, and $\hat{z} = (\hat{t}, \hat{x}', \hat{x}_d) \in \bR \times \overline{\bR^d_+}$. Assume that  \eqref{con-elli}, \eqref{a-b.zero}, and \eqref{a-dd.cond} are satisfied on $((\hat{x}_d -2)^+, \hat{x}_d + 2)$.  If $u \in \sW^{1,2}_2(Q_{2}^+(\hat{z}), d\mu)$ is a strong solution of \eqref{int-Q4-constant-L-0}, then for every $0< \rho < R \leq 2$,
\[
\begin{split}
& \int_{Q_\rho^+(\hat{z})} \big(|Du(z)|^2 + \lambda |u(z)|^2\big) \mu(dz) \leq N(d,\nu, \rho, R) \int_{Q_R^+(\hat{z})}|u(z)|^2 \mu(dz),\\
& \int_{Q_\rho^+(\hat{z})} |u_t(z)|^2 \mu(dz) \leq N(d,\nu, \rho, R) \int_{Q_R^+(\hat{z})} \big(|Du(z)|^2 + \lambda |u(z)|^2\big) \mu(dz).
\end{split}
\]
Moreover, for any $j \in \mathbb{N} \cup \{0\}$, we also have
\[
\begin{split}
& \int_{Q_\rho^+(\hat{z})} |\partial_t^{j+1}u(z)|^2 \mu(dz)+ \int_{Q_\rho^+(\hat{z})} |DD_{x'}\partial_t^{j}u(z)|^2 \mu(dz) \\
&  \leq N(d,\nu, j, \rho, R) \int_{Q_R^+(\hat{z})} \big(|Du(z)|^2 + \lambda |u(z)|^2\big) \mu(dz).
\end{split}
\]
\end{lemma}
\begin{proof} As the equation in \eqref{int-Q4-constant-L-0} can be written in divergence form as in \eqref{div-L-0}, the lemma can be proved by using the standard energy estimates. See, for example,  the proof of \cite[Proposition 4.2]{Dong-Phan}.
\end{proof}

Our next result is the following local  boundary  weighted $L_\infty$  and Lipschitz estimates of solutions.

\begin{lemma}
            \label{Lip-est-lemma}
Let $\nu \in (0, 1]$, $\lambda \geq 0$, and $\alpha <1$ and assume that  \eqref{con-elli}, \eqref{a-b.zero}, and \eqref{a-dd.cond} are satisfied on $(0, 2)$.
If $u \in \sW^{1,2}_2(Q_{2}^+(\hat{z}), d\mu)$ is a strong solution of \eqref{int-Q4-constant-L-0} with $\hat{z} \in \bR \times \overline{\bR^d_+}$, then we have
\begin{equation*}
                                    \begin{split}
& \sup_{z\in Q^+_1(\hat{z})}|x_d^{\alpha}\max\{x_d^{-1}, 1\} u(z)| \leq N \Big(\fint_{Q_{2}^+(\hat{z})} |x_d^\alpha u(z)|^2\mu_1(dz) \Big)^{1/2}, \\
&  \sup_{z\in Q^+_1(\hat{z})}|x_d^\alpha D u(z)| \leq N \Big( \fint_{Q_{2}^+(\hat{z})}\big(|x_d^\alpha Du(z)|^2 +\lambda|x_d^\alpha u(z)|^2\big) \mu_1(dz)\Big)^{1/2},
\end{split}
\end{equation*}
where $N = N(d, \alpha, \nu)>0$.
%\end{itemize}
\end{lemma}
\begin{proof} As already noted, the equation in \eqref{int-Q4-constant-L-0} can be written in the divergence form as in \eqref{div-L-0}. Therefore, Lemma \ref{Lip-est-lemma} follows by applying \cite[Propositions 4.1 and 4.2]{Dong-Phan} to the equation \eqref{div-L-0}.
\end{proof}

We now derive local boundary $L_\infty$-estimates for higher-order derivatives of solutions to the homogeneous equations.
\begin{lemma}
            \label{Higher-Lip-est-lemma}
Let $q\in (1,\infty)$ and $q_1\in [1,\infty)$. Under the assumptions of Lemma \ref{Lip-est-lemma}, if $u \in \sW^{1,2}_{q}(Q_{2}^+(\hat{z}),  x_d^{\alpha q} d\mu_1)$ is a strong solution of \eqref{int-Q4-constant-L-0} and $\hat{z}  =(\hat{z}',0) \in \bR^{d} \times\{0\}$, then for any $j, k \in \mathbb{N} \cup\{0\}$,
\begin{align} \nonumber
& \sup_{z \in Q_1^+(\hat{z})} \big[|x_d^{\alpha-1}D_{x'}^k \partial_t^{j+1} u (z)|  +  |x_d^\alpha D D_{x'}^k \partial_t^j u(z)| + |x_d^{\alpha-1}D_{x'}^k  \partial_t^{j}  u(z)| \big]  \\  \label{H-B-L-infty-u}
&  \leq  N\Big(\fint_{Q_2^+(\hat{z})} |x_d^\alpha D_{x'}^k \partial_t^{j} u(z)|^{q_1} \mu_1(dz) \Big)^{1/q_1},
\end{align}
\begin{align} \nonumber
&  \sup_{z \in Q_1^+(\hat{z})} \big[|\partial_t(x_d^\alpha  D D_{x'}^ku (z))| + |D (x_d^\alpha D D_{x'}^k u(z))|  \big] \\ \label{H-B-L-infty-Du}
&  \leq  N\Big(\fint_{Q_2^+(\hat{z})} \big(x_d^\alpha | DD_{x'}^k u(z)| + \sqrt{\lambda} |D_{x'}^ku(z)|\big)^{q_1}\mu_1(dz) \Big)^{1/q_1},
\end{align}
and
\begin{align}
                                \label{eq4.39}
\|\bM^{\alpha-1} D_{x'}^k u\|_{C^{1/4,1/2}(Q_1^+(\hat{z}))}
\le N\Big(\fint_{Q_2^+(\hat{z})} |x_d^{\alpha-1} D_{x'}^k u(z)|^{q_1} \mu_1(dz) \Big)^{1/q_1}
\end{align}
for $N = N(d, \nu, \alpha, j, k)$. A similar assertion also holds for $\hat{z} = (\hat{z}', \hat{x}_d)$ with $\hat{x}_d > 2$.
\end{lemma}
\begin{proof}  We only prove the boundary estimates since the proof of the interior estimates is simpler. By H\"older's inequality for $q_1>2$ and a standard iteration argument for $q_1\in [1,2)$ (see, for instance, \cite[p. 75]{FHL}), we only need to consider the case when $q_1=2$. By shifting the coordinates, we may also assume that $\hat z=(0,0)$.

We first impose the additional condition that $u \in \sW^{1,2}_2(Q_{2}^+,  d\mu)$.
By using standard argument of finite-difference quotients, we see that $D^{k}_{x'}\partial^{j}_tu$ is still a solution of \eqref{int-Q4-constant-L-0} for $j, k \in \mathbb{N} \cup\{0\}$. Therefore, without loss of generality, we may assume that $j=k=0$. Applying Lemmas \ref{Lip-est-lemma}  and \ref{cappio},  we get
\begin{equation*} %\label{0128.eqn1}
\begin{split}
& \sup_{z\in Q_{1}^+}\Big[ | x_d^{\alpha-1} u_t(z)| +  | x_d^\alpha D u(z)|  + | x_d^{\alpha-1}u(z)|\Big] \\
&  \leq N\Big(\fint_{Q_2^+} |x_d^\alpha u(z)|^{2} \mu_1(dz) \Big)^{1/2},
\end{split}
\end{equation*}
which gives \eqref{H-B-L-infty-u}.

To show \eqref{H-B-L-infty-Du}, as before we may assume that $k=0$. Applying Lemma \ref{Lip-est-lemma} to $u_t$ and then Lemma \ref{cappio}, we get
\begin{align}
                \label{eq4.49}
&\sup_{z \in Q_1^+} |x_d^\alpha  D u_t (z)|\notag\\
&\le  N \bigg( \fint_{Q_{4/3}}\big(|x_d^\alpha Du_t(z)|^2 +\lambda|x_d^\alpha u_t(z)|^2\big) \mu_1(dz)\bigg)^{1/2}\notag\\
&\le  N \bigg( \fint_{Q_{5/3}}|x_d^\alpha u_t(z)|^2 \mu_1(dz)\bigg)^{1/2}\notag\\
&\le  N \bigg( \fint_{Q_{2}}\big(|x_d^\alpha Du(z)|^2 +\lambda|x_d^\alpha u(z)|^2\big) \mu_1(dz)\bigg)^{1/2}.
\end{align}
Applying Lemma \ref{Lip-est-lemma} to $D_{x'}u$ and Lemma \ref{cappio}, we have
\begin{align}
                \label{eq4.57}
&\sup_{z \in Q_1^+} |x_d^\alpha  DD_{x'} u (z)|\notag\\
&\le  N \bigg( \fint_{Q_{3/2}}\big(|x_d^\alpha DD_{x'}u(z)|^2 +\lambda|x_d^\alpha D_{x'}u(z)|^2\big) \mu_1(dz)\bigg)^{1/2}\notag\\
&\le N \bigg( \fint_{Q_{2}}|x_d^\alpha D_{x'}u(z)|^2 \mu_1(dz)\bigg)^{1/2}.
\end{align}
Similarly, applying Lemma \ref{Lip-est-lemma} to $u_t$ and $u$ and then Lemma \ref{cappio}, we have
\begin{align}
                \label{eq5.02}
&\sup_{z \in Q_1^+} |x_d^\alpha  u_t (z)|+\lambda |x_d^\alpha  u (z)|\notag\\
&\le  N \bigg( \fint_{Q_{3/2}}\big(|x_d^\alpha u_t(z)|^2 +\lambda^2|x_d^\alpha u(z)|^2\big) \mu_1(dz)\bigg)^{1/2}\notag\\
&\le N \bigg( \fint_{Q_{2}}\big(|x_d^\alpha Du(z)|^2 +\lambda|x_d^\alpha u(z)|^2\big) \mu_1(dz)\bigg)^{1/2}.
\end{align}
Now we bound $D_d(x_d^\alpha D_d u)$ by using the PDE in \eqref{int-Q4-constant-L-0} and combine \eqref{eq4.49}, \eqref{eq4.57}, and \eqref{eq5.02} to get \eqref{H-B-L-infty-Du}.

Next we prove \eqref{eq4.39}. Again we may assume that $k=0$. In view of \eqref{H-B-L-infty-u}, it suffices to show that for any $(t_1,x_1',x_{1d}), (t_1,x_1',x_{2d})\in Q_1^+$ satisfying $x_{2d}<x_{1d}$, we have
\begin{align}
                    \label{eq4.51}
&|x_{1d}^{\alpha-1}u(t_1,x_1',x_{1d})-x_{2d}^{\alpha-1}u(t_1,x_1',x_{2d})|\notag\\
&\le N(x_{1d}-x_{2d})^{1/2}\Big(\fint_{Q_2^+} |x_d^{\alpha-1} u(z)|^{q_1} \mu_1(dz) \Big)^{1/q_1}.
\end{align}
When $x_{1d}-x_{2d}\le x_{1d}^2/2$, it follows from the mean value formula and \eqref{H-B-L-infty-u} that
\begin{align*}
           %         \label{eq4.51b}
&|x_{1d}^{\alpha-1}u(t_1,x_1',x_{1d})-x_{2d}^{\alpha-1}u(t_1,x_1',x_{2d})|\\
&\le N(x_{1d}-x_{2d})x_{1d}^{-1}\Big(\fint_{Q_2^+} |x_d^{\alpha} u(z)|^{2} \mu_1(dz) \Big)^{1/2}\\
&\le N(x_{1d}-x_{2d})^{1/2}\Big(\fint_{Q_2^+} |x_d^{\alpha-1} u(z)|^{2} \mu_1(dz) \Big)^{1/2},
\end{align*}
which yields \eqref{eq4.51}.
Next we consider the case when $x_{1d}-x_{2d}> x_{1d}^2/2$. By \eqref{H-B-L-infty-Du} and Lemma \ref{cappio}, we have
$$
\sup_{z \in Q_1^+} |D_d (x_d^\alpha D_d u(z))|\le  N\Big(\fint_{Q_2^+} |x_d^{\alpha} u(z)|^{2} \mu_1(dz) \Big)^{1/2}.
$$
Therefore, there exists a bounded function $f=f(z')$ such that for any $z \in Q_1^+$,
$$
|x_d^\alpha D_d u(z)-f(z')|\le Nx_d\Big(\fint_{Q_2^+} |x_d^{\alpha} u(z)|^{2} \mu_1(dz) \Big)^{1/2},
$$
which implies that
$$
| D_d u(z)-x_d^{-\alpha} f(z')|\le Nx_d^{1-\alpha}\Big(\fint_{Q_2^+} |x_d^{\alpha} u(z)|^{2} \mu_1(dz) \Big)^{1/2}.
$$
Using the zero boundary condition, we obtain
$$
| u(z)-(1-\alpha)^{-1}x_d^{1-\alpha} f(z')|\le Nx_d^{2-\alpha}\Big(\fint_{Q_2^+} |x_d^{\alpha} u(z)|^{2} \mu_1(dz) \Big)^{1/2},
$$
which is equivalent to
$$
|x_d^{\alpha-1} u(z)-(1-\alpha)^{-1} f(z')|\le Nx_d\Big(\fint_{Q_2^+} |x_d^{\alpha} u(z)|^{2} \mu_1(dz) \Big)^{1/2}.
$$
Then by the triangle inequality,
\begin{align*}
           %         \label{eq4.51b}
&|x_{1d}^{\alpha-1}u(t_1,x_1',x_{1d})-x_{2d}^{\alpha-1}u(t_1,x_1',x_{2d})|\\
&\le Nx_{1d}\Big(\fint_{Q_2^+} |x_d^{\alpha} u(z)|^{2} \mu_1(dz) \Big)^{1/2}\\
&\le N(x_{1d}-x_{2d})^{1/2}\Big(\fint_{Q_2^+} |x_d^{\alpha-1} u(z)|^{2} \mu_1(dz) \Big)^{1/2},
\end{align*}
which gives \eqref{eq4.51}.

Finally, we remove the additional condition. Observe that if $q \in [2,\infty)$, then by H\"{o}lder's inequality, $u \in \sW^{1,2}_2(Q_{2}^+, d\mu)$. On the other hand, if $q \in (1,2)$, as $u_t$ and $D_{x'}u$ satisfy the same equation as $u$, by using \cite[Corrollary 2.3]{Dong-Phan} for weak solutions to equations in divergence form as in \eqref{div-L-0}, we see that $u_t, D_{x'}u, Du, u  \in L_2(Q_{R}^+, d\mu)$ for any $R<2$. This and Lemma \ref{D-dd-control} imply that $u \in \sW^{1,2}_2(Q_{R}^+, d\mu)$. The lemma is proved.
\end{proof}

We now prove the following result regarding pointwise estimates of solutions to \eqref{int-Q4-constant-L-0} for a more general class of solutions.
\begin{proposition}
            \label{Higher-Lip-est-lemma-s}
Let $p\in (1,\infty)$ and $\beta\in (p-1,(2-\alpha)p-1)$. If
$$
u \in \sW^{1,2}_{p}(Q_{2}^+(\hat{z}), \, x_d^{\alpha p+\beta}\,dz)
$$
is a strong solution of \eqref{int-Q4-constant-L-0} and $\hat{z} =(\hat{z}',0)$, then for any $j, k \in \mathbb{N} \cup\{0\}$,
\begin{align} \nonumber
& \sup_{z \in Q_1^+(\hat{z})} \big[|x_d^{\alpha-1} D_{x'}^k \partial_t^{j+1} u (z)|  +  |x_d^\alpha D D_{x'}^k \partial_t^j u(z)| + |x_d^{\alpha-1}D_{x'}^k  \partial_t^{j}  u(z)| \big]  \\  \label{H-B-L-infty-ubb}
&  \leq  N\fint_{Q_2^+(\hat{z})} |x_d^\alpha D_{x'}^k \partial_t^{j} u(z)| x_d^\beta\, dz,
\end{align}
\begin{align} \nonumber
&  \sup_{z \in Q_1^+(\hat{z})} \big[|\partial_t(x_d^\alpha  D D_{x'}^ku (z))| + |D (x_d^\alpha D D_{x'}^k u(z))|  \big] \\ \label{H-B-L-infty-Dubb}
&  \leq  N\fint_{Q_2^+(\hat{z})} \big(x_d^\alpha  |DD_{x'}^k u(z)| + \sqrt{\lambda} |D_{x'}^ku(z)|\big)x_d^\beta\, dz,
\end{align}
and
\begin{align}
                                \label{eq4.39bb}
\|\bM^{\alpha-1} D_{x'}^k u\|_{C^{1/4,1/2}(Q_1^+(\hat{z}))}
\le N\fint_{Q_2^+(\hat{z})} |x_d^{\alpha-1} D_{x'}^k u(z)|x_d^\beta dz
\end{align}
for $N = N(d, \nu, \alpha, j, k,\beta)$. A similar assertion also holds for $\hat{z} = (\hat{z}', \hat{x}_d)$ with $\hat{x}_d > 2$.
\end{proposition}
\begin{proof}  As before, we only consider the boundary case. Without loss of generality, we may assume that $\hat z=0$. The proposition follows directly from Lemma \ref{Higher-Lip-est-lemma} if $\beta \le -\alpha$. Next we consider the case when $\beta>-\alpha$. We first impose the additional condition that $u \in \sW^{1,2}_{p}(Q_{3/2}^+, \, x_d^{\alpha p}d\mu_1)$. Then  applying Lemma \ref{Higher-Lip-est-lemma}, we obtain
\begin{align*}
& \sup_{z \in Q_1^+} \big[|x_d^{\alpha} D_{x'}^k \partial_t^{j} u (z)|+
|x_d^{\alpha-1} D_{x'}^k \partial_t^{j+1} u (z)|  +  |x_d^\alpha D D_{x'}^k \partial_t^j u(z)| + |x_d^{\alpha-1}D_{x'}^k  \partial_t^{j}  u(z)| \big]  \\
&  \leq  N\fint_{Q_{4/3}^+} |x_d^\alpha D_{x'}^k \partial_t^{j} u(z)| x_d^{-\alpha}\, dz\\
&  \leq  \varepsilon \fint_{Q_{4/3}^+} |x_d^\alpha D_{x'}^k \partial_t^{j} u(z)| x_d^{\beta_0}\,dz
+N\varepsilon^{(\beta+ \alpha)/(\alpha +\beta_0)} \fint_{Q_{4/3}^+} |x_d^\alpha D_{x'}^k \partial_t^{j} u(z)| x_d^{\beta}\, dz\\
&  \leq  N\varepsilon \sup_{Q_{4/3}^+} |x_d^\alpha D_{x'}^k \partial_t^{j} u(z)|+N\varepsilon^{(\beta+\alpha)/(\alpha +\beta_0)} \fint_{Q_{4/3}^+} |x_d^\alpha D_{x'}^k \partial_t^{j} u(z)| x_d^{\beta}\, dz,
\end{align*}
where $\beta_0\in (-1,-\alpha)$, and we also used Young's inequality in the  second inequality. By an iteration argument (see, for instance, \cite[Lemma 4.3]{FHL}), we obtain \eqref{H-B-L-infty-ubb}. The estimates \eqref{H-B-L-infty-Dubb}  and \eqref{eq4.39bb} can be proved similarly.

Next, we remove the additional condition that $u \in \sW^{1,2}_{p}(Q_{3/2}^+ ,x_d^{\alpha p}\, d\mu_1)$. By taking the standard mollification with respect to $t$ and $x'$, and then taking the limit, without loss of generality, we may assume that  $u$ is smooth in $x'$ and $t$. We claim that
\begin{equation}
                            \label{eq5.42}
u \in \sW^{1,2}_{p}(Q_{3/2}^+, \, x_d^{\alpha p+\beta'}dz)\quad \forall\, \beta'>-1,
\end{equation}
which implies \eqref{H-B-L-infty-ubb}, \eqref{H-B-L-infty-Dubb}, and \eqref{eq4.39bb}
in the general case. To this end,
we take a smooth cutoff function $\eta\in C_0^\infty((-2,2))$ satisfying
$ \eta\equiv 1$ in $(-5/3,5/3)$.
Then by applying the weighted Hardy's inequality (Lemma \ref{w-hardy-inq} (ii)) to $v:=\eta(x_d) x_d^\alpha D_d u$, we get
\begin{align} \nonumber
& \|\bM^{\alpha-1}D_d u\|_{L_p(Q_2'\times (0,5/3),x_d^\beta)} \le \|\bM^{-1}v\|_{L_p(Q_2'\times \bR^+,x_d^\beta)} \\  \label{eq5.34}
&\le N\|D_d v\|_{L_p(Q_2'\times \bR^+,x_d^\beta)}<\infty,
\end{align}
where we used the condition $\beta>p-1$. By applying the weighted inequality Lemma \ref{w-hardy-inq} (i) and using the zero boundary condition, we have
$$
\|\bM^{\alpha-2} u\|_{L_p(Q_2'\times (0,5/3),x_d^\beta)}
\le N\|\bM^{\alpha-1}D_d u\|_{L_p(Q_2'\times (0,5/3),x_d^\beta)}<\infty,
$$
where we used the condition $\beta<(2-\alpha)p-1$.
Therefore,
$$
u\in L_p(Q^+_{5/3},x_d^{(\alpha-2) p+\beta}\,dz).
$$
By applying the same argument to $D_{x'} u$, we also have
$$
D_{x'}u\in L_p(Q^+_{5/3},x_d^{(\alpha-2) p+\beta}\,dz).
$$
This together with \eqref{eq5.34} imply that
$$
\bM^{-1}D_{x'}u, \quad Du\in L_p(Q^+_{5/3},x_d^{\alpha p+\beta_1}\,dz),\quad\text{where}\,\,\beta_1=\beta-p.
$$
Again as $u_t$ and $D_{x'} u$ satisfy the same conditions as $u$, we obtain
$$
u_t, DD_{x'}u\in L_p(Q^+_{5/3},x_d^{\alpha p+\beta_1}\,dz).
$$
Then from \eqref{int-Q4-constant-L-0} and Lemma \ref{D-dd-control}, we get $D(\bM^\alpha D_du)\in L_p(Q^+_{5/3},x_d^{\beta_1}\,dz)$, and thus $u\in \sW^{1,2}_{p}(Q_{3/2}^+, \, x_d^{\alpha p+\beta_1}dz)$.
If $\beta_1>p-1$, we can repeat this procedure to get $u\in \sW^{1,2}_{p}(Q_{r_k}^+, \, x_d^{\alpha p+\beta_k}dz)$, where $\beta_k=\beta-kp$ and $\{r_k\}\subset (3/2,2)$ is a finite sequence of decreasing numbers. Let $k_0$ be the integer such that $\beta_{k_0-1}>p-1$ and $\beta_{k_0}\le p-1$. We then have that  $u\in \sW^{1,2}_{p}(Q_{r_{k_0}}^+, \, x_d^{\alpha p+\hat \beta}dz)$ for any $\hat\beta>p-1$. By repeating this procedure once again, we prove the claim \eqref{eq5.42}. The proposition is proved.
\end{proof}

%We now conclude the section with the following remark.

\subsection{Mixed-norm \texorpdfstring{$L_p$}{Lp}-estimates for non-homogeneous equations} \label{non-simple-coeffients}

The main result of the section is following theorem on the existence and estimate of solutions in $\hat\sW^{1,2}_{q,p}(\Omega_T,  x_d^{\gamma}\, dz)$ to \eqref{constant-L-0} with $\gamma \in (\alpha p -1, 2p-1\big)$.
\begin{theorem} \label{thrm-s6} Let $\nu \in (0, 1]$, $p, q \in (1,\infty)$, $\alpha \in (-\infty,1)$, and $\gamma \in \big(\alpha p -1, 2p-1\big)$ be constants. Assume that $\overline{a}_{ij}$ satisfies  \eqref{con-elli} and \eqref{a-dd.cond}, and $\overline{a}_0, \overline{c}$ satisfy \eqref{a-b.zero}. Then, for any $f \in L_{q,p}(\Omega_T, x_d^{\gamma}\, dz)$ and $\lambda > 0$, there exists a unique strong solution $u \in\hat{\sW}^{1,2}_{q,p}(\Omega_T,  x_d^{\gamma}\, dz)$ to \eqref{constant-L-0} and $u$ satisfies
\begin{equation} \label{all-s-constant-L-p}
\begin{split}
&\|u_t \|
+\|D_{x} D_{x'} u\| + \|D_{d}^2 u + \alpha \bM^{-1} D_{d} u\|\\
& +\|\bM^{-1} D_{x'} u\| +\lambda^{1/2}\|Du\|
+\lambda\|u\| +\lambda^{1/2}\|\bM^{-1} u\| \le N\|f\|,
\end{split}
\end{equation}
where $\|\cdot\| = \|\cdot\|_{L_{q,p}(\Omega_T, x_d^\gamma dz)}$, and $N= N(d, \nu,\alpha, p,q, \gamma) >0$.
\end{theorem}

Observe that Theorem \ref{thm0} is a special case of Theorem \ref{thrm-s6}. Before proving the theorem, let us recall some notation of functional spaces used in \cite{Dong-Phan}. Let $p \in [1, \infty)$, $S, T \in [-\infty,+\infty]$ with $S<T$, $\cD \subset \bR^d_+$ be open, and $\tau \in \bR$. We define
\[
\begin{split}
& \bH_{p}^{-1}( (S,T)\times \cD, x_d^\tau\, dz) \\
& =\big\{u:\, u  =  D_iF_i+\bM^{-1}F_0 +f\ \ \text{for some}\ f\in L_{p}( (S,T)\times \cD, x_d^\tau\, dz)\\
& \qquad \ F= (F_0,\ldots,F_d) \in L_{p}((S,T)\times \cD, x_d^\tau\, dz)^{d +1}\big\},
\end{split}
\]
which is equipped with the norm
\begin{align*}
\|u\|_{\bH_{p}^{-1}((S,T)\times \cD, x_d^\tau\, dz)} &=\inf\big\{\|F\|_{L_{p}((S,T)\times \cD, x_d^\tau\, dz)}
+\|f\|_{L_{p}((S,T)\times \cD,x_d^\tau\, dz)}: \\
& \quad \qquad \ u= D_iF_i+\bM^{-1}F_0 +f\big\}.
\end{align*}
We also define the space
\[
\begin{split}
& \sH_{p}^1((S,T)\times \cD, x_d^\tau\, dz)\\
& =\big\{u \in L_p((S, T),  \sW^1_p(\cD, x_d^\tau\, dx)):    u_t\in  \bH_{p}^{-1}( (S,T)\times \cD, x_d^\tau\, dz)\big\}
\end{split}
\]
equipped with the norm
\begin{align*}
&\|u\|_{\sH_{p}^1((S,T)\times \cD, x_d^\tau\, dz)} \\
&= \|u\|_{L_{p}((S,T)\times \cD,x_d^\tau\, dz)}
+ \|Du\|_{L_{p}((S,T)\times \cD,x_d^\tau\, dz)}+\|u_t\|_{\bH_{p}^{-1}((S,T)\times \cD, x_d^\tau\, dz)}.
\end{align*}

Next, we state and prove the following lemma.
\begin{lemma}
\label{s-W-2-p.constant.eqn} Let $\nu \in (0, 1]$, $p \in (1,\infty)$, $\alpha \in (-\infty,1)$, and $\gamma \in \big((\alpha +1) p -1, 2p-1\big)$ be constants. Assume that $\overline{a}_{ij}$ satisfies  \eqref{con-elli} and \eqref{a-dd.cond}, and $\overline{a}_0, \overline{c}$ satisfy \eqref{a-b.zero}. Then, for any $f \in L_p(\Omega_T, x_d^{\gamma}\, dz)$ and $\lambda > 0$, there exists a unique strong solution $u \in  \sH^1_p(\Omega_T,  x_d^{\gamma-p}\, dz)\cap \hat{\sW}^{1,2}_p(\Omega_T,  x_d^{\gamma}\, dz)$ to \eqref{constant-L-0}, which satisfies
\begin{equation} \label{s-constant-L-p}
\begin{split}
&\|u_t \|
+\|D^2 u\|
+\lambda^{1/2}\|Du\|
 +\|\bM^{-1}Du\| \\
& +\lambda\|u\|  +\sqrt{\lambda}\|\bM^{-1}u\| + \|\bM^{-2}u\| \le N\|f\|,
\end{split}
\end{equation}
where $\|\cdot\| = \|\cdot\|_{L_p(\Omega_T, x_d^\gamma dz)}$ and $N= N(d, \nu,\alpha, p,\gamma) >0$.
\end{lemma}
\begin{proof}
By applying the scaling argument mentioned in Remark \ref{remark2}, we may assume $\lambda =1$. Also, by using a density argument, we assume that $f$ is compactly supported in $\Omega_T$. Let $\gamma_1 = \gamma - (\alpha +1) p$, and
\begin{equation*}
              %                      \label{eq11.36}
F(t,x)=-\int_{x_d}^\infty s^\alpha f(t,x',s)\,ds.
\end{equation*}
Note that $D_d F(z)= x_d^\alpha f(z)$ for $z=(t,x) \in \Omega_T$ with $x = (x', x_d) \in \bR^d_+$.  By the weighted Hardy's inequality (Lemma \ref{w-hardy-inq} (ii)), we have
\begin{align}
            \label{eq4.39cc}
&\|F\|_{L_{p}(\Omega_T,x_d^{\gamma_1}dz)}
=\|\bM^{-1}F\|_{L_{p}(\Omega_T,x_d^{\gamma_1+p}dz)}\notag\\
&\le N\|D_d F\|_{L_{p}(\Omega_T,x_d^{\gamma_1+p}dz)}
\le N\|\bM^\alpha f\|_{L_{p}(\Omega_T,x_d^{\gamma_1+p}dz)} = \| f\|_{L_{p}(\Omega_T,x_d^{\gamma}dz)},
\end{align}
where we used the condition that $\gamma_1+p=\gamma-\alpha p>p-1$. By the assumptions on the coefficients $\overline{a}_{ij}$, we can write the operator $\mathcal{L}_0$ in divergence form as in \eqref{div-L-0}. Then, the equation \eqref{constant-L-0} becomes
\begin{equation} \label{div-L-0-1}
x_d^{\alpha}\big(\overline{a}_0(x_d)u_t + \overline{c}(x_d) u \big)- D_i[x_d^\alpha \overline{a}_{ij}(x_d) D_j u + F_i]   =  0  \quad \text{in} \quad \Omega_T,
\end{equation}
where $F_i = 0$ for $i =1, 2,\ldots, d-1$ and $F_d =F$. Note also that by the assumption on $\gamma$, we see that $\gamma_1\in (-1,(1-\alpha) p-1)$. Then, by applying \cite[Theorem 2.5 and Remark 2.6]{Dong-Phan} to the equation \eqref{div-L-0-1}, we see that there is a unique weak solution $u\in \sH^1_p(\Omega_T,  x_d^{\alpha p+\gamma_1}\,dz)$ of \eqref{div-L-0-1}, and
$$
\|\bM^\alpha Du\|_{L_{p}(\Omega_T,x_d^{\gamma_1}dz)} + \|\bM^\alpha u\|_{L_{p}(\Omega_T,x_d^{\gamma_1}dz)} \le N\|F\|_{L_{p}(\Omega_T,x_d^{\gamma_1}dz)},
$$
which together with \eqref{eq4.39cc} implies that
\begin{equation} \label{Du-0226}
\|Du\|_{L_{p}(\Omega_T,x_d^{\gamma -p}dz)} + \|u\|_{L_{p}(\Omega_T,x_d^{\gamma - p}dz)} \le N\|f\|_{L_{p}(\Omega_T,x_d^{\gamma}dz)}.
\end{equation}
As $\gamma - p+1 < p$, by using the weighted Hardy's inequality (Lemma \ref{w-hardy-inq} (i)) and the boundary condition $u =0$ on $\{x_d =0\}$, we get
\begin{align} \notag
 \|u \|_{L_{p}(\Omega_T,x_d^{\gamma-2p}dz)} &= \|\bM^{-1}u\|_{L_{p}(\Omega_T,x_d^{\gamma -p}dz)}  \le N\|D_du\|_{L_{p}(\Omega_T,x_d^{\gamma -p}dz)} \\\label{u-op-wei-0226}
& \le N\|f\|_{L_{p}(\Omega_T,x_d^{\gamma}dz)}.
\end{align}

Now we prove the estimate \eqref{s-constant-L-p} and conclude that $u\in \hat{\sW}^{1,2}_p(\Omega_T,  x_d^{\gamma}\, dz)$. We follow an idea introduced in \cite[Lemma 2.2]{Kr99} (see also \cite[Theorem 3.5]{DK15}).  Let $\zeta \in C_0^\infty(\bR_+)$ be non-negative such that
\[
\int_0^\infty |\zeta (s)|^p s^{-\gamma-1} ds =1, \quad
\int_0^\infty |\zeta'(s)|^p s^{p-\gamma-1} ds =N_1 <\infty,
\]
and
\[
\int_0^\infty |\zeta''(s)|^p s^{2p-\gamma  -1} ds = N_2 <\infty.
\]
For each $r >0$, let $\zeta_r(s) =\zeta(rs)$, where $s \in \bR_+$. We note that with a suitable assumption on the integrability of a given function $v: \Omega_T \rightarrow \bR$, using the Fubini theorem, we have
\begin{equation} \label{weight-kry}
\begin{split}
& \int_0^\infty \left(\int_{\Omega_T}|\zeta_r(x_d) v(z)|^p\, dz\right) r^{-\gamma-1}\, dr =\int_{\Omega_T} |v(z)|^p x_d^{\gamma}\, dz, \\
& \int_0^\infty \left(\int_{\Omega_T}|\zeta_r'(x_d) v(z)|^p\,  dz\right) r^{-\gamma-1}\, dr =N_1\int_{\Omega_T} |v(z)|^p x_d^{\gamma -p}\, dz, \\
& \int_0^\infty \left(\int_{\Omega_T}|\zeta_r''(x_d) v(z)|^p\,  dz\right) r^{-\gamma-1}\, dr =N_2\int_{\Omega_T} |v(z)|^p x_d^{\gamma-2p}\, dz.
\end{split}
\end{equation}
Next, for each fixed $r>0$, let $w(t,x) =\zeta_r(x_d) u(t,x)$, where $(t, x) \in \Omega_T$ with $x= (x', x_d) \in \bR^{d-1} \times \bR_+$. Then $w$ solves the following equation with uniformly elliptic coefficients
\begin{equation} \label{w-eqn-0307}
\overline{a}_0w_t + \overline{c} w - \overline{a}_{ij}D_{ij}w = g \quad \text{in} \quad \Omega_T\quad \text{and} \quad w =0 \quad \text{on} \quad \{x_d =0\},
\end{equation}
where
\[
\begin{split}
g(t,x) & =\zeta_r(x_d) f(t,x) + \frac{\alpha}{x_d} \big(\zeta_r(x_d) D_du + \zeta'_r(x_d) u  \big) \\
& \quad - \sum_{i=1}^{d-1} \overline{a}_{id} \zeta'_r (x_d)D_i u  - 2 \zeta_r' D_d u - \zeta''_r(x_d) u.
\end{split}
\]
We first prove \eqref{s-constant-L-p}  with the assumption that $u \in  \sW^{1,2}_p(\Omega_T, x_d^{\alpha p +\gamma_1}dz)$.  Under this extra assumption and as $\zeta_r$ is compactly supported in $(0, \infty)$, we see that $w \in W^{1,2}_p(\Omega_T)$, where $W^{1,2}_p(\Omega_T)$ is the usual parabolic Sobolev space. Then by applying the  $W^{1, 2}_{p}$-estimate for \eqref{w-eqn-0307} (see, for instance, \cite{D12}),  we obtain
\[
\|w\|_{W^{1,2}_p(\Omega_T)}  \leq N \| g\|_{L_p(\Omega_T)}.
\]
From this, the definition of $g$, and a simple manipulation, we obtain
\[
\begin{split}
& \|\zeta_r u \|_{L_p(\Omega_T)} +  \|\zeta_r Du\|_{L_p(\Omega_T)}  + \|\zeta_r D^2u\|_{L_p(\Omega_T)} +\|\zeta_r u_t \|_{L_p(\Omega_T)} \\
& \leq N\Big[ \| \zeta f\|_{L_p(\Omega_T)} +  \|\zeta''_r u\|_{L_p(\Omega_T)}  + \|\zeta_r' Du\|_{L_p(\Omega_T)} + \|\zeta_r' u\|_{L_p(\Omega_T)} \\
& \qquad + \|\zeta_r Du\|_{L_p(\Omega_T, x_d^{-p}dz)} + \|\zeta'_r u\|_{L_p(\Omega_T, x_d^{-p}dz)} \Big].
\end{split}
\]
Then, by raising this last estimate to the power $p$, multiplying it with $r^{-\gamma-1}$, and then integrating with respect to $r$ on $(0,\infty)$, we obtain
\[
\begin{split}
& \| u\|_{L_p(\Omega_T, x_d^{\gamma}dz)} + \|Du \|_{L_p(\Omega_T, x_d^{\gamma}dz)} + \| D^2u\|_{L_p(\Omega_T, x_d^{\gamma}dz)}  + \| u_t\|_{L_p(\Omega_T, x_d^{\gamma}dz)}\\
 & \leq N \Big[ \|f\|_{L_p(\Omega_T, x_d^{\gamma}dz)} + \|Du\|_{L_p(\Omega_T, x_d^{\gamma -p}dz)}  \\
 & \qquad + \|u\|_{L_p(\Omega_T, x_d^{\gamma -p}dz)}  + \|u\|_{L_p(\Omega_T, x_d^{\gamma -2p}dz)}\Big].
 \end{split}
\]
where we also used \eqref{weight-kry}. It then follows from the last estimate, \eqref{Du-0226}, and \eqref{u-op-wei-0226} that
\[
\begin{split}
&  \| u\|_{L_p(\Omega_T, x_d^{\gamma}dz)} + \|Du \|_{L_p(\Omega_T, x_d^{\gamma}dz)} + \| D^2u\|_{L_p(\Omega_T, x_d^{\gamma}dz)}  + \| u_t\|_{L_p(\Omega_T, x_d^{\gamma}dz)}\\
&  \quad + \|Du\|_{L_p(\Omega_T, x_d^{\gamma -p}dz)}  + \|u\|_{L_p(\Omega_T, x_d^{\gamma -p}dz)} + \|u\|_{L_p(\Omega_T, x_d^{\gamma -2p}dz)} \\
& \leq N  \|f\|_{L_p(\Omega_T, x_d^{\gamma}dz)}.
\end{split}
\]

It remains to remove the extra assumption that $u \in \sW^{1,2}_p(\Omega_T, x_d^{\alpha p + \gamma_1}dz)$. By mollifying the equation \eqref{constant-L-0} in $t$ and $x'$ and applying \cite[Theorem 2.5 and Remark 2.6]{Dong-Phan} to the equations of $u^{(\varepsilon)}_t$ and $D_{x'}u^{(\varepsilon)}$, we obtain
$$
u^{(\varepsilon)}, u_t^{(\varepsilon)}, D_{x'}u^{(\varepsilon)}, DD_{x'}u^{(\varepsilon)} \in L_p(\Omega_T, x_d^{\alpha p+\gamma_1} dz).
$$
This and Lemma \ref{D-dd-control} imply that $u^{(\varepsilon)} \in \sW^{1,2}_p(\Omega_T, x_d^{\alpha p+ \gamma_1}dz)$ is a strong solution of \eqref{constant-L-0} with $f^{(\varepsilon)}$ in place of $f$. From this, we apply the a priori estimate \eqref{s-constant-L-p} that we just proved for $u^{(\varepsilon)}$ and pass to the limit as $\varepsilon \rightarrow 0^+$ to obtain the estimate \eqref{s-constant-L-p} for $u$. The proof of the theorem is completed.
\end{proof}

\begin{remark} \label{rm-47} Note that in Lemma \ref{s-W-2-p.constant.eqn} we estimated $D_d^2 u$ because $\gamma$ is large. The lower bound $(\alpha+1)p-1$ of $\gamma$ is sharp by considering the example $u(x) = x^{1-\alpha}$ and $u''(x) = -\alpha (1-\alpha) x^{-1-\alpha}$.  When $\alpha=0$, this result is also consistent with the previous results in \cite{K99,DK15}.
\end{remark}

From Lemma \ref{s-W-2-p.constant.eqn} and Proposition \ref{Higher-Lip-est-lemma-s}, we derive the following mean oscillation estimate.
\begin{corollary}
            \label{s-osc-lemma-2022}
Let $\nu \in (0, 1]$, $p_0 \in (1, \infty)$, $\alpha \in (-\infty,1)$, and $\gamma_0 \in (p_0-1, (2-\alpha)p_0-1)$ be constants. Let $\lambda>0$, $\rho>0$,  and $\hat{z} = (\hat{t}, \hat{x}', \hat{x}_d) \in \overline{\Omega}_T$. Assume that  \eqref{con-elli}, \eqref{a-b.zero}, and \eqref{a-dd.cond} are satisfied. If $f \in L_{p_0}(Q_{8\rho}^+(\hat{z}), x_d^{\alpha p_0}d\mu_2)$ and $u \in \hat{\sW}^{1,2}_{p_0}(Q_{8\rho}^+(\hat{z}), x_d^{\alpha p_0}d\mu_2)$ is a strong solution to the equation
 \[
 \left\{
 \begin{array}{cccl}
  \mathcal{L}_0 u  & = & f &  \quad \text{in} \quad Q_{6\rho}^+(\hat{z}), \\
  u & = & 0 &  \quad \text{on} \quad  Q_{6\rho}(\hat{z})\cap \{x_d=0\} \quad \text{if} \quad \hat{x}_d \leq 6\rho,
  \end{array} \right.
 \]
where $\mu_2(s) = s^{\gamma_0}$ for $s \in \bR_+$, then
\begin{equation}
                    \label{eq10.36}
\begin{split}
&  \fint_{Q_{\kappa\rho}^+(\hat{z})} |U- (U)_{Q_{\kappa \rho}^+(\hat{z})}| \,\mu_2(dz)\\
 & \leq N\kappa^{1/2}  (|U|)_{Q_{8\rho}^+(\hat{z})} + N\kappa^{-(d +2+\gamma_0)/p_0}\left( |\bM^\alpha f|^{p_0} \right)_{Q_{8\rho}^+(\hat{z})}^{1/p_0}
 \end{split}
\end{equation}
for any $\kappa \in (0,1)$, where
\begin{equation}
                            \label{eq5.50}
U = \bM^{\alpha}(u_t, D D_{x'} u, \bM^{-1}D_{x'}u, \sqrt{\lambda} Du, \lambda u,\sqrt\lambda \bM^{-1} u),
\end{equation}
$(\cdot)_{Q}$ is defined in \eqref{k-meanb}, and $N=N(\nu, d, \alpha, p_0, \gamma_0)>0$.
\end{corollary}

\begin{proof}
By the assumption on $\gamma_0$, we see that $\overline{\gamma}: = \alpha p_0 + \gamma_0 \in ( (1+\alpha) p_0 - 1, 2p_0 -1)$. Then, by applying Lemma \ref{s-W-2-p.constant.eqn} with $\overline{\gamma}$ and $p_0$ in place of $\gamma$ and $p$, there is a solution $v \in \hat{\sW}^{1,2}_{p_0}(\Omega_T, x_d^{\overline{\gamma}}\, dz)$ to
\[
\left\{
\begin{array} {cccl}
\mathcal{L}_0 v & = & f {\bf 1}_{Q^+_{8\rho}(\hat z)}& \quad \text{in} \quad \Omega_T \\
 v & = & 0 & \quad \text{on} \quad \{x_d =0\}
\end{array} \right.
\]
satisfying
\begin{equation} \label{v-est-0214}
\begin{split}
& \norm{v_t}_{L_{p_0}(\Omega_T, x_d^{\overline\gamma}\, dz)} +  \norm{D^2 v}_{L_{p_0}(\Omega_T, x_d^{\overline{\gamma}}\, dz)} + \norm{\bM^{-1} Dv}_{L_{p_0}(\Omega_T, x_d^{\overline\gamma}\, dz)} \\
&\quad + \sqrt{\lambda} \norm{Dv}_{L_{p_0}(\Omega_T, x_d^{\overline\gamma}\, dz)} + \lambda \norm{v}_{L_{p_0}(\Omega_T, x_d^{\overline\gamma}\, dz)}
 +\sqrt\lambda \norm{\bM^{-1} v}_{L_{p_0}(\Omega_T, x_d^{\overline\gamma}\, dz)}\\
&\leq N  \norm{f}_{L_{p_0}(Q_{8\rho}^+(\hat{z}), x_d^{\overline\gamma}\, dz)}.
\end{split}
\end{equation}
From the definition of $\mu_2$ in \eqref{mu-def} and the definition of $(\cdot)_{Q}$ in \eqref{k-meanb}, we  see that \eqref{v-est-0214}  implies
\begin{equation} \label{V-0309-est}
\begin{split}
&\big (|V|^{p_0}\big)_{Q_{\kappa \rho}^+(\hat{z})}^{1/p_0} \leq N \kappa^{-(d+2 +\gamma_0)/p_0} \big(|\bM^\alpha f|^{p_0}\big)_{Q_{8\rho}^+(\hat{z})}^{1/p_0}, \\
& \big(|V|^{p_0}\big)_{Q_{8 \rho}^+(\hat{z})}^{1/p_0} \leq N \big(|\bM^\alpha f|^{p_0}\big)_{Q_{8\rho}^+(\hat{z})}^{1/p_0},
\end{split}
\end{equation}
where
$$
V = \bM^{\alpha} (v_t, DD_{x'}v,  \bM^{-1}D_{x'}v, \sqrt{\lambda}Dv, \lambda v, \sqrt\lambda \bM^{-1}v).
$$
Now, let $w = u-v \in \hat{\sW}^{1,2}_{p_0}(Q_{8\rho}^+(\hat{z}), x_d^{\bar{\gamma}}\, dz)$, and note that $w$ solves the equation
\[
\mathcal{L}_0 w = 0 \quad \text{in} \quad Q_{8\rho}^+(\hat{z})
\]
and $w =0$ on $\{x_d =0\} \cap Q_{8\rho}(\hat{z})$. Then, we use the mean value theorem  and Proposition \ref{Higher-Lip-est-lemma-s} for $w$ with $ \gamma_0$ in place of $\beta$ and $p_0$ in place of $p$ to obtain
\[
\fint_{Q_{\kappa \rho}^+(\hat{z})} |W - \big(W\big)_{Q_{\kappa \rho}^+(\hat{z})}|\, \mu_2(dz) \leq N \kappa^{1/2} \big(W\big)_{Q_{8 \rho}^+(\hat{z})},
\]
where
$$
W = \bM^{\alpha}(w_t, DD_{x'}w,  \bM^{-1} D_{x'}w, \sqrt{\lambda}Dw, \lambda w,  \sqrt\lambda \bM^{-1}w).
$$
From this last estimate, \eqref{V-0309-est}, and the triangle inequality, we obtain the desired estimate. The proof of the lemma is completed.
\end{proof}
\smallskip
Finally, we provide the proof of Theorem \ref{thrm-s6}.
\begin{proof}[Proof of Theorem \ref{thrm-s6}]
We first prove the a priori estimate \eqref{all-s-constant-L-p} for a given $u \in \hat{\sW}^{1,2}_{q,p}(\Omega_T,  x_d^{\gamma}\, dz)$ solving \eqref{constant-L-0}.  Let  $p_0 \in \big(1, \min\{p, q \}\big)$ be sufficiently closed to $1$,  and  $\gamma_0 \in (p_0 -1, (2-\alpha) p_0 -1)$ be sufficiently close to $(2-\alpha) p_0 -1$ such that
\begin{equation*}
                %            \label{eq1.04-1}
\gamma': =\gamma - \alpha p -\gamma_0 \in \big(-1-\gamma_0,(1+\gamma_0)(p/p_{0}-1)\big),
\end{equation*}
which implies that
$$
x_d^{\gamma'} \in  A_{p/p_0}(\bR_+, \mu_2)
\subset A_{p}(\bR_+, \mu_2).
$$
Applying Corollary \ref{s-osc-lemma-2022}, we see that for every $\rho>0$ and $\hat{z} \in \Omega_T$, we have
\begin{equation*}
\begin{split}
&  \fint_{Q_{\kappa\rho}^+(\hat{z})} |U- (U)_{Q_{\kappa \rho}^+(\hat{z})}| \,\mu_2(dz)\\
 & \leq N\kappa^{1/2}  (|U|)_{Q_{8\rho}^+(\hat{z})} + N\kappa^{-(d +2+\gamma_0)/p_0}\left( |\bM^\alpha f|^{p_0} \right)_{Q_{8\rho}^+(\hat{z})}^{1/p_0},
 \end{split}
\end{equation*}
where $\kappa \in (0,1)$, $U$  is defined in \eqref{eq5.50}, and $N=N(\nu, d, \alpha, p_0, \gamma_0)>0$. It then follows that
\[
U^{\#} \leq N \kappa^{1/2} \M(|U|) + N\kappa^{-(d +2+\gamma_0)/p_0} \M ( |\bM^\alpha f|^{p_0} )^{1/p_0} \quad \text{on} \quad \Omega_T,
\]
where the sharp function and the maximal function are defined in \eqref{sharp-def} and \eqref{maxi-def}, respectively.
Then, by using the weighted Fefferman-Stein theorem and the Hardy-Littlewood theorem for $\M$ with the weight $x_d^{\gamma'}$ (Theorem \ref{FS-thm}), we obtain the estimate
\[
\begin{split}
\|U\|_{L_{q,p}(\Omega_T,  x_d^{\gamma'}d\mu_2)} & \le N\Big[\kappa^{1/2} \|U\|_{L_{q,p}(\Omega_T,  x_d^{\gamma'}d\mu_2)} \\
& \qquad + \kappa^{-(d +2+\gamma_0)/p_0}  \| \bM^\alpha f\|_{L_{q,p}(\Omega_T, x_d^{\gamma'}d\mu_2)}\Big].
\end{split}
\]
From this, and by choosing $\kappa>0$  sufficiently small, we obtain
\[
\|U\|_{L_{q,p}(\Omega_T,  x_d^{\gamma'}d\mu_2)} \le N \|\bM^\alpha f\|_{L_{q,p}(\Omega_T,  x_d^{\gamma'}d\mu_2)}.
\]
Then, by using Lemma \ref{D-dd-control} for \eqref{constant-L-0}, we obtain the estimate \eqref{all-s-constant-L-p}.

We now prove the solvability of \eqref{constant-L-0} in $\hat{\sW}^{1,2}_{q,p}(\Omega_T,  x_d^{\gamma}\, dz)$. We split the proof into two cases.
\\
\noindent
{\bf Case I}:  $\gamma \in (\alpha p-1, p-1)$.  We write $\mathcal{L}_0$ in divergence form so that the equation \eqref{constant-L-0} becomes
\begin{equation} \label{div-L-0-2}
x_d^{\alpha}\big(\overline{a}_0(x_d)u_t + \lambda \overline{c}(x_d) u \big)- D_i[x_d^\alpha \overline{a}_{ij}(x_d) D_j u]   = x_d^\alpha f  \quad \text{in} \quad \Omega_T.
\end{equation}
It follows from \cite[Theorem 2.5 and Remark 2.6]{Dong-Phan} that there is a unique weak solution $u \in \sH^1_{q,p}(\Omega_T, x_d^{\gamma}dz)$ of \eqref{div-L-0-2}. Then, as in the proof Lemma \ref{s-W-2-p.constant.eqn}, by mollifying the equation \eqref{div-L-0-2} in $(t,x')$ and using Lemma \ref{D-dd-control}, we see that $u^{(\varepsilon)} \in \sW^{1,2}_{q,p}(\Omega_T, x_d^\gamma dz)$ is a strong solution of \eqref{constant-L-0} with  $f^{(\varepsilon)}$ in place of $f$.  Moreover, since $\gamma<p-1$, by using the weighted Hardy's inequality (Lemma \ref{w-hardy-inq} (i)), we know that
$$
x_d^{-1}u^{(\varepsilon)} ,x_d^{-1}D_{x'}u^{(\varepsilon)} \in  L_{q,p}(\Omega_T, x_d^\gamma dz).
$$
From this, we apply the a priori estimate \eqref{all-s-constant-L-p} that we just proved for $u^{\varepsilon}$ and then pass to the limit as $\varepsilon \rightarrow 0^+$, we see that $u \in \hat{\sW}^{1,2}_{q,p}(\Omega_T, x_d^\gamma dz)$ and it solves \eqref{constant-L-0}.
\ \\
\noindent
{\bf Case II}: $\gamma \in [ p-1, 2p -1)$. For each $f \in L_{q, p}(\Omega_T, x_d^\gamma dz)$, let $\{f_k\} \subset C_{0}^\infty(\Omega_T)$ such that
\[
\lim_{k\rightarrow \infty}\|f_k - f\|_{L_{q, p}(\Omega_T, x_d^\gamma dz)} =0.
\]
We take a $\gamma_1\in (\alpha p-1 ,p-1)$ so that $\gamma_1<\gamma$. Since $f_k$ has compact support in $\Omega_T$ and smooth, we have $f_k \in L_{q, p}(\Omega_T, x_d^{\gamma_{1}} dz)$. Therefore, by {\bf Case I}, there exists a solution $u_k \in \hat{\sW}^{1,2}_{q,p}(\Omega_T, x_d^{\gamma_1} dz)$ of \eqref{constant-L-0} with  $f_k$ in place of $f$.

Now, if the sequence $\{u_k\}\subset \hat{\sW}^{1,2}_{q,p}(\Omega_T, x_d^{\gamma} dz)$, then by the a priori estimate \eqref{all-s-constant-L-p}, we obtain
\[
\begin{split}
&\|\partial_t u_k \|
+\|DD_{x'} u_k\| + \|D_{d}^2 u + \alpha \bM^{-1} D_{d} u_k\| +\|\bM^{-1}D_{x'} u_k\|\\
&\quad +\lambda^{1/2}\|Du_k\|
+\lambda\|u_{k}\| +\lambda^{1/2}\|\bM^{-1} u_k\| \le N\|f_k\|,
\end{split}
\]
where $\|\cdot\| = \|\cdot\|_{L_{q,p}(\Omega_T, x_d^\gamma dz)}$. It follows from the linearity of \eqref{constant-L-0} and the convergence of $\{f_k\}$ in $L_{q, p}(\Omega_T, x_d^\gamma dz)$  that $\{u_k\}$ is Cauchy in $\hat{\sW}^{1,2}_{q,p}(\Omega_T, x_d^{\gamma} dz)$.  Let $u \in \hat{\sW}^{1,2}_{q,p}(\Omega_T, x_d^{\gamma} dz)$ be the limit of  $\{u_k\}$ in $\hat{\sW}^{1,2}_{q,p}(\Omega_T, x_d^{\gamma} dz)$. Then, by passing to the limit, we see that $u$ solves \eqref{constant-L-0}.

It remains to prove that $u_k \in \hat{\sW}^{1,2}_{q,p}(\Omega_T, x_d^{\gamma} dz)$ for all $k$. Let $k$ be fixed, and let $R_0 >0$ be sufficiently large such that
\[
\text{supp}(f_k) \subset D_{R_0} := (-\infty,  T) \times \bR^{d-1} \times (0,R_0).
\]
As $\gamma_1 < \gamma$ and $u_k \in \hat{\sW}^{1,2}_{q,p}(\Omega_T, x_d^{\gamma_1} dz)$,  we see that $u_k \in \hat{\sW}^{1,2}_{q,p}(D_{2R_0}, x_d^{\gamma} dz)$. Hence, we only need to show that
\[
\|u_k\|_{\hat{\sW}^{1,2}_{q,p}(\Omega_T \setminus D_{R_0}, x_d^{\gamma} dz)} <\infty.
\]
For each $l \in \mathbb{N} \cup \{0\}$, let $\eta_l=\eta_l(x_d)$ be a smooth function such that $\eta_l \equiv 0$ in $(-2^l R_0,2^l R_0)$, $\eta_l \equiv 1$ outside $(-2^{l+1} R_0,2^{l+1} R_0)$, and
\[
\|D^{k} \eta_l \|_{L_\infty} \leq N_0 2^{-kl}, \quad k=0,1,2, \quad \forall l \geq 0.
\]
Let $w_l = u_k \eta_l$, which is a solution to the equation
\[\mathcal{L}_0 w_l = g \quad \text{in} \quad \Omega_T \quad \text{and} \quad w = 0 \quad \text{on} \quad \{x_d =0\},\]
where
\[
g = - \overline{a}_{ij}(D_i u \delta_{jd} \eta_l' + D_j u \delta_{id} \eta_l'+u  \delta_{id}\delta_{jd}\eta''_l)- \alpha x_d^{-1} u \eta_l'
\]
Therefore, by the a priori estimate \eqref{all-s-constant-L-p}, we have
\begin{equation*}% \label{est.0313}
\begin{split}
&\|\partial_t w_l \|_{L_{q, p}(\Omega_T, x_d^{\gamma_1} dz)}
+\|D D_{x'} w_l\|_{L_{q, p}(\Omega_T, x_d^{\gamma_1} dz)} \\
& + \|D_{d}^2 w_l + \alpha \bM^{-1} D_{d} w_l\|_{L_{q, p}(\Omega_T, x_d^{\gamma_1} dz)}  +\lambda^{1/2}\|Dw_l\|_{L_{q, p}(\Omega_T, x_d^{\gamma_1} dz)} \\
& + \|\bM^{-1} D_{x'}w_l\|_{L_{q, p}(\Omega_T, x_d^{\gamma_1} dz)} + \lambda^{1/2} \|\bM^{-1} w_l\|_{L_{q, p}(\Omega_T, x_d^{\gamma_1} dz)} \\
& +\lambda \|w_l \|_{L_{q, p}(\Omega_T, x_d^{\gamma_1} dz)}    \le N\|g\|_{L_{q, p}(\Omega_T, x_d^{\gamma_1} dz)}.
\end{split}
\end{equation*}
Observe that
\[
\begin{split}
 \|g\|_{L_{q, p}(\Omega_T, x_d^{\gamma_1} dz)}
 & \leq  N 2^{-l}\Big(\|Du\|_{L_{q,p}(D_{2^{l+1}R_0} \setminus D_{2^{l}R_0}, x_d^{\gamma_1} dz)} \\
 & \qquad + 2^{-l}\|u  \|_{L_{q,p}(D_{2^{l+1}R_0} \setminus D_{2^{l}R_0}, x_d^{\gamma_1} dz)} \Big).
\end{split}
\]
Combining the two inequalities above, we get
\[
\|u\|_{\hat{\sW}^{1,2}_{q, p}(D_{2^{l+2}R_0} \setminus D_{2^{l+1}R_0}, x_d^{\gamma_1} dz)} \leq N 2^{-l} \|u\|_{\hat{\sW}^{1,2}_{q, p}(D_{2^{l+1}R_0} \setminus D_{2^lR_0}, x_d^{\gamma_1} dz)},
\]
where $N>0$ depends also on $\lambda$, but is independent of $l$. By iterating this estimate, we obtain
\[
\|u\|_{\hat{\sW}^{1,2}_{q, p}(D_{2^{l+1}R_0} \setminus D_{2^{l}R_0}, x_d^{\gamma_1} dz)} \leq N^l 2^{-l(l-1)/2} \|u\|_{\hat{\sW}^{1,2}_{q, p}(D_{2R_0}, x_d^{\gamma_1} dz)}, \quad \forall \ l \geq 0.
\]
This implies
\[
\begin{split}
\|u\|_{\hat{\sW}^{1,2}_{q, p}(\Omega_T \setminus D_{R_0}, x_d^\gamma dz)} & = \sum_{l=0}^\infty \|u\|_{\hat{\sW}^{1,2}_{q, p}(D_{2^{l+1}R_0} \setminus D_{2^{l}R_0}, x_d^\gamma dz)} \\
& \leq  \|u\|_{\hat{\sW}^{1,2}_{q, p}(D_{2R_0}, x_d^{\gamma_1} dz)} \sum_{l=0}^\infty N^{l} 2^{-l(l-1)/2} (2^l R_0)^{(\gamma-\gamma_1)/p} \\
& \leq  N\|u\|_{\hat{\sW}^{1,2}_{q, p}(D_{R_0}, x_d^{\gamma_1} dz)} < \infty.
\end{split}
\]
The proof is completed.
\end{proof}

We now conclude this section with the following lemma, which improves Corollary \ref{s-osc-lemma-2022} as the lower bound of $\gamma_0$ does not depend on $p_0$. The lemma is used in the next section.

\begin{lemma} \label{s-osc-0228}
Corollary \ref{s-osc-lemma-2022} still holds when $\gamma_0 \in (-1, (2-\alpha)p_0+1)$ provided that the term $\kappa^{-(d+2+\gamma_0)/p_0}$ in \eqref{eq10.36} is replaced with  $\kappa^{-(d+2+\gamma_0^+)/p_0}$, where $\gamma_0^+=\max\{\gamma_0, 0\}$.
\end{lemma}
\begin{proof} We repeat the proof of Corollary \ref{s-osc-lemma-2022} by using Theorem \ref{thrm-s6} instead of Lemma \ref{s-W-2-p.constant.eqn}.
\end{proof}

\begin{remark} \label{a-corr-result} As the discussion right after \eqref{a-dd.cond}, the results in this section hold when \eqref{a-dd.cond} is replaced with \eqref{eq4.31}.
\end{remark}
%==============
\section{Equations with partially VMO coefficients} \label{VMO-section}

This section is to prove Theorems \ref{para-main.theorem}, \ref{elli-main.theorem}, and Corollary \ref{cor2.3}. We shall first study the equation \eqref{main-eqn} which is a parabolic equation in non-divergence form with singular coefficients:
\begin{equation} \label{eqn.variable-coeff}
\left\{
\begin{array}{cccl}
 \mathcal{L} u  &  = & f & \quad \text{in} \quad \Omega_T, \\
u & = & 0 & \quad \text{on} \quad \{x_d =0\},
\end{array}  \right.
\end{equation}
where $\cL$ is defined in \eqref{L.def}.
Note that we can always divide both sides of the PDE in \eqref{eqn.variable-coeff} by $a_{dd}$ and replace $\nu$ in \eqref{ellipticity} and \eqref{a-b.zero0} with $\nu^2$. Therefore, it is convenient to assume that
\begin{equation} \label{extension-type-matrix}
a_{dd}=1.
\end{equation}
We first state and prove a lemma about the oscillation estimate  for  solutions to non-homogeneous equations.
\begin{lemma} \label{lemma.3.3} Let $\nu \in (0,1)$, $p_0 \in (1, \infty)$, $\alpha \in (-\infty,1), \gamma_0 \in (-1, (2-\alpha)p_0 -1)$, $p \in (p_0, \infty)$ and assume that \eqref{ellipticity}, \eqref{a-b.zero0}, and \eqref{extension-type-matrix} are satisfied. Let $\lambda>0$ and $\rho, \rho_1, \rho_0 \in (0,1)$,  $\hat{z} = (\hat{t}, \hat{x}', \hat{x}_d) \in \overline{\Omega}_T$, $t_1 \in \bR$ and $f \in L_{p_0}(Q_{8\rho}^+(\hat{z}), x_d^{p_0\alpha} d\mu_2)$. Assume that $u \in  \hat{\sW}^{1,2}_{p}(Q_{8\rho}^+(\hat{z}), x_d^{p\alpha}d\mu_2)$ vanishing outside $(t_1 -(\rho_0\rho_1)^2, t_1]$ is a strong solution to the equation
 \[
 \left\{
 \begin{array}{cccl}
 \mathcal{L} u  & = &f & \quad \text{in} \quad Q_{6\rho}^+(\hat{z}),\\
u  & = & 0 & \quad \text{on} \quad  Q_{6\rho}(\hat{z})\cap\{x_d=0\}  \quad \text{if} \quad \hat{x}_d \leq 6\rho.
 \end{array} \right.
 \]
 Then, for any $\kappa \in (0,1)$, it holds that
\begin{equation}
                    \label{eq8.43}
\begin{split}
& \fint_{Q_{\kappa \rho}^+(\hat{z})} | U - (U)_{Q_{\kappa \rho}^+(\hat{z})} | \,\mu_2(dz)\\
&\leq N\kappa^{1/2} \left(|U|\right)_{Q_{8\rho}^+(\hat{z})} +  N\kappa^{-(d +2+\gamma_0^+)} \rho_1^{2(1-1/p_0)}\left( | U |^{p_0} \right)_{Q_{8\rho}^+(\hat{z})}^{1/p_0}  \\
&  \quad + N\kappa^{-\frac{d +2+\gamma_0^+}{p_0}}\Big[\left( |\bM^\alpha f|^{p_0} \right)_{Q_{8\rho}^+(\hat{z})}^{1/p_0}+a^{\#}_{\rho_0}(\hat{z})^{\frac{1}{p_0} -\frac{1}{p}} \left(|U|^p  \right)_{Q_{8\rho}^+(\hat{z})}^{1/p} \Big],
\end{split}
\end{equation}
where $U$ is defined in \eqref{eq5.50} and $N=N(d, \nu, p, p_0, \alpha, \gamma_0) >0$.
\end{lemma}
\begin{proof}
We discuss two cases depending on $8\rho > \rho_0$ or $8\rho \leq \rho_0$. \\
{\bf Case I}:  $8\rho > \rho_0$.  By using the doubling property of $\mu_2$ and H\"older's inequality, we have
\[
\begin{split}
 & \fint_{Q_{\kappa \rho}^+(\hat{z})} | U - (U)_{Q_{\kappa \rho}^+(\hat{z})} | \,\mu_2(dz)\leq   N \kappa^{-(d +2+\gamma_0^+)}(| U|)_{Q_{8\rho}^+(\hat{z})}  \\
 & \leq N \kappa^{-(d +2+\gamma_0^+)} \left(  \mathbf{1}_{(t_1 -(\rho_0\rho_1)^2, t_1]}\right)_{Q_{8 \rho}^+(\hat{z})}^{1-{1}/{p_0}} \left(  | U |^{p_0} \right)_{Q_{8\rho}^+(\hat{z})}^{{1}/{p_0}} \\
 &\leq N \kappa^{-(d +2+\gamma_0^+)} \rho_1^{2(1-1/p_0)}\left(  | U |^{p_0}\right)_{Q_{8\rho}^+(\hat{z})}^{1/p_0},
 \end{split}
\]
where $N = N(d, p_0, \gamma_0)>0$.

\noindent
{\bf Case 2}: $8\rho \leq \rho_0$. Recall that  $[a_0]_{8\rho, \hat{z}}(x_d)$, $[c]_{8\rho, \hat{z}}(x_d)$,  $[a_{ij}]_{8\rho, \hat{z}}(x_d)$, and $(a_{dj})_{8\rho, \hat{z}}$ are defined as in \eqref{average-a}-\eqref{k-mean}, where $i, j \in \{1, 2,\ldots, d\}$.
Denote
\[
\mathcal{L}_{\rho,\hat{z}} u = [a_0]_{8\rho, \hat{z}}u_t + \lambda  [c]_{8\rho, \hat{z}}u -    \overline{a}_{ij} (x_d) D_{ij} u - \frac{\alpha}{x_d} \overline{a}_{dj} {D_{j}} u,
\]
where
\begin{equation}
                \label{eq8.41}
\overline{a}_{ij} (x_d)= \left \{
\begin{array}{ll}
[a_{ij}]_{8\rho, \hat{z}}(x_d), & \quad \text{for} \quad i=1,2,\ldots, d-1, \quad j =1, 2,\ldots, d \\
(a_{dj})_{8\rho, \hat{z}}, & \quad \text{for} \quad i=d, \quad j =1, 2,\ldots, d
\end{array} \right.
\end{equation}
and
\[
\begin{split}
F_1(z)  & = \sum_{i=1}^{d-1}\sum_{j=1}^d(a_{ij} - [a_{ij}]_{8\rho, \hat{z}}) D_{ij} u(z), \\
F_2(z) & = \sum_{j=1}^{d-1}(a_{dj} - (a_{dj})_{8\rho, \hat{z}}) (D_{dj} u(z) +\alpha x_d^{-1}D_j u(z)), \\
F_3(z)  & =  ([a_{0}]_{8\rho, \hat{z}}-a_{0}) u_t(z)
+\lambda ([c]_{8\rho, \hat{z}}-c) u(z).
\end{split}
\]
Under the condition \eqref{extension-type-matrix},  $u$ satisfies
\begin{equation*}
 \mathcal{L}_{\rho, \hat{z}} u(t,x) =  f(t,x) + \sum_{i=1}^3 F_i(t,x)  \quad \text{in} \quad Q_{6\rho}^+(\hat{z})
\end{equation*}
with the boundary condition $u =0$ on $\{x_d =0\}$ if $\hat{x}_d \leq 6\rho$. Also, as $1 < p_0 < p$ and $\gamma_0>-1$, by using H\"{o}lder's inequality, $u \in \hat{\sW}^{1,2}_{p_0}(Q_{8\rho}^+(\hat{z}), x_d^{p_0\alpha}d\mu_2)$. From \eqref{extension-type-matrix} and the definitions in \eqref{eq8.41},  the coefficient matrix $(\overline{a}_{ij})$ satisfies \eqref{eq4.31}. Hence, as explained in Remark \ref{a-corr-result}, we can apply Lemma \ref{s-osc-0228} and infer that
\begin{equation} \label{F-oss}
\begin{split}
& \fint_{Q_{\kappa\rho}^+(\hat{z})} |U- (U)_{Q_{\kappa \rho}^+(\hat{z})}| \,\mu_2(dz)\\
& \leq N\Big[\kappa^{1/2}  (|U|)_{Q_{8\rho}^+(\hat{z})}+ \kappa^{-(d +2+\gamma_0^+)/{p_0}}\left(  |\bM^\alpha f|^{p_0} \right)_{Q_{8\rho}^+(\hat{z})}^{1/p_0}  \\
& \qquad +   \kappa^{-(d +2+\gamma_0^+)/p_0}\sum_{i=1}^{3}\left(  |\bM^\alpha F_i|^{p_0}\right)_{Q_{8\rho}^+(\hat{z})}^{1/p_0}\Big],
\end{split}
\end{equation}
where $N = N( d,\nu, \alpha, \gamma_0, p_0) >0$. We now bound the last term on the right-hand side of \eqref{F-oss}.  By H\"{o}lder's inequality and the boundedness of $(a_{ij})$ in \eqref{ellipticity}, the definition of $a_\rho^{\#}$ in \eqref{a-sharp}, and \eqref{extension-type-matrix},
\[
\begin{split}
 &\left(|\bM^\alpha F_1|^{p_0}\right)_{Q_{8\rho}^+(\hat{z})}^{1/p_0}\\
 &\leq  \Big(  |a_{ij}(z) -[a_{ij}]_{8\rho, \hat{z}}(x_d) |^{pp_0/(p-p_0)} \Big)_{Q_{8\rho}^+(\hat{z})}^{1/p_0-1/p} \left( |\bM^\alpha  D D_{x'}u|^{p} \right)_{Q_{8\rho}^+(\hat{z})}^{1/p} \\
& \leq  N \left( |a_{ij}(z) -[a_{ij}]_{8\rho, \hat{z}}(x_d) |\right)_{Q_{8\rho}^+(\hat{z})}^{1/p_0-1/p} \left( |\bM^\alpha  DD_{x'}u|^p \right)_{Q_{8\rho}^+(\hat{z})}^{1/p} \\
& = N a^{\#}_{\rho_0}(\hat z)^{1/p_0-1/p}\left(  |\bM^\alpha  DD_{x'}u|^p \right)_{Q_{8\rho}^+(\hat{z})}^{1/p},
\end{split}
\]
where $N = N(d, \nu, p, p_0)>0$. Similarly, we also have
\[
\begin{split}
&\left( |\bM^\alpha  F_2|^{p_0} \right)_{Q_{8\rho}^+(\hat{z})}^{1/p_0} \leq  N a^{\#}_{\rho_0}(\hat z)^{1/p_0-1/p}\left(  |\bM^\alpha  DD_{x'}u|^p  +|\bM^{\alpha-1}  D_{x'}u|^p
\right)_{Q_{8\rho}^+(\hat{z})}^{1/p},\\
& \left( |\bM^\alpha  F_3|^{p_0} \right)_{Q_{8\rho}^+(\hat{z})}^{1/p_0} \leq N a^{\#}_{\rho_0}(\hat{z})^{1/p_0-1/p}\left( |\bM^\alpha  u_t|^p+\lambda^p |\bM^\alpha  u|^p  \right)_{Q_{8\rho}^+(\hat{z})}^{1/p}.
\end{split}
\]
By plugging the estimates of $F_k$ for $k=1,2,3$ into \eqref{F-oss}, we obtain \eqref{eq8.43}. From the above two cases, the lemma is proved.
\end{proof}
\begin{proposition}  \label{small-spt.thrm} Let $\nu$, $T$, $p$, $q$, $K$, $\alpha$, $\gamma_0$, $\rho_0$, and $\omega$ be as in Theorem \ref{para-main.theorem}. There exist sufficiently small positive constants $$
\delta = \delta (d, \nu, \alpha, p, q, \gamma_0,  K)\quad \text{and} \quad
\rho_1 = \rho_1(d, \nu, \alpha, p, q, \gamma_0, K)
$$
such that, under the conditions \eqref{ellipticity}, \eqref{a-b.zero0}, and \eqref{para-VMO}, the following statement holds. Let $\lambda >0$ and $f \in L_{q,p}(\Omega_T, x_d^{\alpha p}{\omega}\ d\mu_2)$. If $u \in \hat{\sW}^{1,2}_{q,p}(\Omega,   x_d^{\alpha p}\omega\, d\mu_2)$ vanishes outside $(t_1 - (\rho_0 \rho_1)^2, t_1]$ for some $t_1 \in \bR$ and satisfies \eqref{eqn.variable-coeff}, then
\[
\begin{split}
& \|u_t\|_{L_{q,p}(\Omega_T,  x_d^{\alpha p}\omega\,d\mu_2)} + \|DD_{x'} u\|_{L_{q,p}(\Omega_T,  x_d^{\alpha p} \omega\,d\mu_2)}
 + \|\bM^{-1} D_{x'} u\|_{L_{q,p}(\Omega_T,  x_d^{\alpha p} \omega\,d\mu_2)}\\
& \quad + \|D_d(\bM^\alpha D_d u)\|_{L_{q,p}(\Omega_T,  \omega\,d\mu_2)}
  + \sqrt{\lambda} \|Du\|_{L_{q,p}(\Omega_T,  x_d^{\alpha p}\omega\,d\mu_2)} \\
&\quad
  + \lambda \|u\|_{L_{q,p}(\Omega_T, x_d^{\alpha p}\omega\,d\mu_2)}
+ \sqrt\lambda \|\bM^{-1} u\|_{L_{q,p}(\Omega_T, x_d^{\alpha p}\omega\,d\mu_2)} \\
&  \leq N(d, \nu, \alpha, p, q, \gamma_0, K) \|f\|_{L_{q,p}(\Omega_T, x_d^{\alpha p}\omega\,d\mu_2)} .
\end{split}
\]
\end{proposition}

\begin{proof} Without loss of generality, we may assume that \eqref{extension-type-matrix} holds. As $\omega_0 \in A_q((-\infty,T))$ and $\omega_1  \in A_p(\bR^d_+, d\mu_2)$, by  the reverse H\"older's inequality \cite[Theorem 3.2]{MS1981}, we find $p_1=p_1(d,p,q,\gamma_0,K)\in (1,\min(p,q))$ such that
\begin{equation}
							\label{eq0605_13}
\omega_0 \in A_{q/p_1}((-\infty,T)),\quad
\omega_1  \in A_{p/p_1}(\bR^d_+, \mu_2).
\end{equation}
Note that because $\gamma_0 \in (-1, 1-\alpha]$, we can choose $p_2 \in (1,p_1)$ such that $\gamma_0 \in (-1, (2-\alpha)p_2-1)$. Applying Lemma \ref{lemma.3.3} with $p_2$  in place of $p_0$ and $p_1$ in place of  $p$,   we have
\[
\begin{split}
U^{\#} \leq & N\Big[\kappa^{1/2} \mathcal{M}(|U|) + \kappa^{-(d +2+\gamma_0^+)} \rho_1^{2(1-1/p_2)}\mathcal{M}(|U|^{p_2})^{1/p_2} \\
&   +\kappa^{-\frac{d +2+\gamma_0^+}{p_2}} \mathcal{M} (|\bM^\alpha f|^{p_2})^{1/p_2} +  \kappa^{-\frac{d +2+\gamma_0^+}{p_2}} \delta^{\frac{1}{p_2} -\frac{1}{p_1}}  \mathcal{M}(|U|^{p_1})^{1/p_1} \Big] \quad \text{in} \quad \overline{\Omega_T}
\end{split}
\]
 for any $\kappa\in (0,1)$, where  $N = N(\nu, d, p_1, p_2, \alpha, \gamma_0) >0$. Therefore, it follows from Theorem \ref{FS-thm} that
\[
\begin{split}
&\norm{U}_{L_{q,p}} \leq  N \Big[ \kappa^{1/2} \| \mathcal{M}(|U|)\|_{L_{q,p}}  + \kappa^{-(d +2+\gamma_0^+)} \rho_1^{2(1-1/p_2)}\| \mathcal{M}(|U|^{p_2})^{1/p_2}\|_{L_{q,p}}   \\
&  \qquad + \kappa^{-\frac{d +2+\gamma_0^+}{p_2}} \| \mathcal{M} (|\bM^\alpha f|^{p_2})^{\frac 1 {p_2}}\|_{L_{q,p}}+  \kappa^{-\frac{d +2+\gamma_0^+}{p_2}} \delta^{\frac{1}{p_2} -\frac{1}{p_1}} \|\mathcal{M}(|U|^{p_1})^{\frac 1 {p_1}}\|_{L_{q,p}}
\Big],
\end{split}
\]
where $N = N(d,\nu, p, q, \alpha, \gamma_0, K)>0$ and $L_{q,p}=L_{q,p}(\Omega_T, \omega\,d\mu_2)$. Then, from \eqref{eq0605_13} and Theorem \ref{FS-thm} again, we get
\[
\begin{split}
&\norm{U}_{L_{q,p}}  \leq N \Big[ \Big(\kappa^{1/2} + \kappa^{-(d +2+\gamma_0^+) }\rho_1^{2(1-1/p_2)}\Big) \| U\|_{L_{q,p}}+ \kappa^{-\frac{d +2+\gamma_0^+}{p_2}} \|\bM^\alpha f\|_{L_{q,p}}  \\
& \quad +  \kappa^{-\frac{d +2+\gamma_0^+}{p_2}} \delta^{\frac{1}{p_2} -\frac{1}{p_1}}  \|U\|_{L_{q,p}} \Big].
\end{split}
\]
Now, by choosing $\kappa$ sufficiently small and then $\delta$ and $\rho_1$ sufficiently small depending on $d,\nu, p, q,\alpha, \gamma_0$, and $K$ such that
\[
N\Big (\kappa^{1/2} + \kappa^{-(d +2+\gamma_0^+)} \rho_1^{2(1-1/p_2)} +  \kappa^{-\frac{d +2+\gamma_0^+}{p_2}} \delta^{\frac{1}{p_2} -\frac{1}{p}}\Big) <1/2,
\]
we obtain
\[
\begin{split}
& \norm{U}_{L_{q,p}}   \leq  N(d, \nu, p, q, \alpha, \gamma_0, K)  \|\bM^\alpha f\|_{L_{q,p}}.
\end{split}
\]
This and Lemma \ref{D-dd-control} prove the assertion of the proposition.
\end{proof}

Now, we are ready to prove Theorem \ref{para-main.theorem}.

\begin{proof}[Proof of Theorem \ref{para-main.theorem}]
We  first prove the estimate \eqref{main-para.est}. Let $u\in \hat{\sW}^{1,2}_{q, p}(\Omega_T, x_d^{p \alpha}\omega\ d\mu_2)$ be a strong solution of \eqref{main-eqn}. We apply Proposition \ref{small-spt.thrm} and a partition of unity argument in the time variable. Let $\xi \in C_0^\infty(\bR)$  be a non-negative standard cut-off function vanishing outside $(-\rho_0^2\rho_1^2, 0]$ and satisfying
\begin{equation} \label{normalized}
\int_{\bR} \xi^q(t)\ dt =1\quad\text{and}\quad  \int_{\bR} (\xi'(t))^q\ dt \le N(\rho_0\rho_1)^{-2q},
\end{equation}
where $\rho_1>0$ is from Proposition \ref{small-spt.thrm}. For a given $s\in \bR$, let $w_s(t,x) = u(t,x)\xi(t-s)$. We see that $w_s$ is a strong solution of
\[
\left\{
\begin{array}{cccl}
 \mathcal{L} w_s & =  & F_s & \quad \text{in} \quad \Omega_T  \\
 w_s(t,x', 0) & = & 0 & \quad \text{for} \quad (t, x') \in (-\infty, T) \times \bR^{d-1},
\end{array}  \right.
\]
where
$$F_s(t,x) = f (t,x) \xi(t-s) +  a_0 u(t,x) \xi_t (t-s).
$$
As $w_s$ vanishes outside $(s-\rho_0^2\rho_1^2, s] \times \bR^d_{+}$, by Proposition \ref{small-spt.thrm}, we have
\begin{equation} \label{w-s.est}
\begin{split}
& \|\partial_t w_s\|_{L_{q,p}} + \sqrt{\lambda} \|Dw_s\|_{L_{q,p}} + \|DD_{x'} w_s\|_{L_{q,p}}
+ \|\bM^{-1} D_{x'} w_s\|_{L_{q,p}} \\
& \quad + \|\bM^{-\alpha}D_d(\bM^\alpha D_{d} w_s)\|_{L_{q,p}} + \lambda \|w_s\|_{L_{q,p}} + \sqrt\lambda \|\bM^{-1} w_s\|_{L_{q,p}} \\
& \leq N \|F_s\|_{L_{q,p}},
\end{split}
\end{equation}
where $N = N (d, \nu, \alpha, \gamma_0, p, q, K)$ and $L_{q,p} = L_{q,p}(\Omega_T, x_d^{\alpha p}\omega\, d\mu_2)$.  From \eqref{normalized}, for any integer $k\ge 0$ and $\tau \in \bR$, we have
\[
\|\bM^\tau D^k_x u\|_{L_{q,p}}^q = \int_{\bR}\|\bM^\tau D_x^k w_s\|^q_{L_{q,p}} \ ds.
\]
Also, it follows from $u_t \xi(t-s) = \partial_t w_s - u \xi_t(t-s)$ that
\[
\|u_t\|^{q}_{L_{q,p}}  \leq  N\int_{\bR} \|\partial_t  w_s\|^q_{L_{q,p}} \ ds + N (\rho_0\rho_1)^{-2 q} \|u\|^{ q}_{L_{q,p}}.
\]
From the last two estimates
and by integrating  the $q$-th power of  \eqref{w-s.est} with respect to $ s$, we conclude that
\[
\begin{split}
& \|u_t\|_{L_{q,p}}  + \sqrt{\lambda} \|D u\|_{L_{q,p}} + \|D D_{x'}u\|_{L_{q,p}}
+ \|\bM^{-1} D_{x'}u\|_{L_{q,p}}\\
&\quad + \|\bM^{-\alpha}D_d(\bM^\alpha D_{d} u)\|_{L_{q,p}}
+ \lambda \|u\|_{L_{q, p} }  + \sqrt\lambda \|\bM^{-1} u\|_{L_{q, p} }\\
& \leq N \|f\|_{L_{q,p}} +  N(\rho_0\rho_1)^{-2} \|u\|_{L_{q,p}},
\end{split}
\]
where $N = N(d, \nu, \alpha, p, q, \gamma_0, K)>0$. Then, by choosing
$$
\lambda_0 = \lambda(d, \nu, p, q, \alpha, \gamma_0, K) = 2  N \rho_1^{-2},
$$
we obtain \eqref{main-para.est} provided that $\lambda \geq \lambda_0 \rho_0^{-2}$.

Observe that the estimate \eqref{main-para.est} also implies the uniqueness of solutions. It then remains to prove the existence of the solution. We split the proof into two steps.
\\
\noindent
{\bf Step I}: Assume $p=q$ and $\omega\equiv1$. We use the method of continuity. Consider the operator
\[
\mathcal{L}_{\tau} u  = (1-\tau)\Big[\partial_t - \Delta -\frac{\alpha}{x_d} D_d + \lambda \Big] u + \tau \mathcal{L} u
\]
with $\tau \in [0,1]$. It is simple to check that the coefficients in $\mathcal{L}_\tau$ satisfy all assumptions in Theorem \ref{para-main.theorem} uniformly in $\tau \in [0,1]$. Then, using the solvability result in Theorem \ref{thrm-s6} and the a priori estimate \eqref{main-para.est} that we just proved, we can apply the method of continuity to obtain the solvability of \eqref{main-eqn} with $\lambda \geq \lambda_0\rho_0^{-2}$, where $\lambda_0 =\lambda_0 (d, \nu, p,q, \alpha, \gamma_0, K)>0$ is defined in the proof of \eqref{main-para.est}. For details, see for example,  \cite[Theorem 1.3.4, p. 15]{Krylov-book} and the proof of \cite[Theorem 1.1]{Dong-Phan-RMI}.
\ \\
\noindent
{\bf Step II}: We consider the general case with $p, q \in (1, \infty)$ and $\omega$ as in the statement of Theorem \ref{para-main.theorem}. We follow the approach in \cite[Section 8]{Dong-Kim-18}. Let $p_1 > \max\{p,q\}$ be sufficiently large and let $\varepsilon_1, \varepsilon_2 \in (0,1)$ be sufficiently small depending on $K, p, q$ and $\gamma_0$ such that
\begin{equation*}
1-\frac{p}{p_1} = \frac{1}{1+\varepsilon_1} \quad \text{and} \quad 1 - \frac{q}{p_1} = \frac{1}{1+\varepsilon_2},
\end{equation*}
and both $\omega_1^{1+\varepsilon_1}$ and $\omega_0^{1+\varepsilon_2}$ are locally integrable and satisfy the doubling property. Precisely, there is $N_0>0$ such that
\begin{equation} \label{omega-0}
\int_{\Gamma_{2r}(t_0)} \omega_0^{1+\varepsilon_2}(s)\, ds \leq N_0 \int_{\Gamma_{r}(t_0)} \omega_0^{1+\varepsilon_2}(s)\, ds
\end{equation}
for any $r>0$ and $t_0 \in \bR$, where $\Gamma_{r}(t_0) = (t_0 -r^2, \min\{t_0 + r^2, T\})$. Similarly
\begin{equation} \label{omega-1-0308}
\int_{B_{2r}^+(x_0)} \omega_1^{1+\varepsilon_1}(x)\, d\mu_2 \leq N_0\int_{B_{r}^+(x_0)} \omega_1^{1+\varepsilon_1}(x)\, d\mu_2
\end{equation}
for any $r >0$ and any $x_0 \in \overline{\bR^d_+}$.

Now, let $\{f_k\}$ be a sequence in $C_0^\infty(\Omega_T)$ such that
\begin{equation} \label{f-k-converge-0227}
\lim_{k\rightarrow \infty} \|f_k - f\|_{L_{q,p}(\Omega_T, x_d^{\alpha p}\omega\, d\mu_2)} =0.
\end{equation}
By Step I, for each $k \in \mathbb{N}$, we can find a solution $u_k \in \hat{\sW}^{1,2}_{p_1}(\Omega_T, x_d^{\alpha p_1}\, d\mu_2)$ of  \eqref{main-eqn} when $f$ is replaced with $f_k$, where $\lambda \geq \lambda_0 \rho_0^{-2}$ for $\lambda_0 = \lambda_0(d, \nu, p_1, p_1,\alpha, \gamma_0, K)>0$. Observe that if the sequence $\{u_k\}$ is in $\hat{\sW}^{1,2}_{q,p}(\Omega_T, x_d^{\alpha p}\omega\, d\mu_2)$, then by applying the a priori estimate \eqref{main-para.est} that we just proved, we have
\[
\|u_k\|_{\hat{\sW}^{1,2}_{q,p}(\Omega_T, x_d^{\alpha p}\omega\, d\mu_2)} \leq N \|f_k\|_{L_{p,q}(\Omega_T, x_d^{\alpha p}\omega\, d\mu_2)}
\]
and
\[
\|u_k - u_l \|_{\hat{\sW}^{1,2}_{q,p}(\Omega_T, x_d^{\alpha p}\omega\, d\mu_2)} \leq N \|f_k - f_l\|_{L_{p,q}(\Omega_T, x_d^{\alpha p}\omega\, d\mu_2)}
\]
for all $k, l \in \mathbb{N}$ and for $N = N(d, p, q, \alpha, \gamma_0, \lambda, K)>0$. Then, by \eqref{f-k-converge-0227}, we see that $\{u_k\}$ is Cauchy in $\hat{\sW}^{1,2}_{q,p}(\Omega_T, x_d^{\alpha p}\omega\, d\mu_2)$. Let $u \in \hat{\sW}^{1,2}_{q,p}(\Omega_T, x_d^{\alpha p}\omega\, d\mu_2)$ be the limit of  $\{u_k\}$ in $\hat{\sW}^{1,2}_{q,p}(\Omega_T, x_d^{\alpha p}\omega\, d\mu_2)$. By passing to the limit as $k \rightarrow \infty$ in the equation of $u_k$, we see that $u \in \hat{\sW}^{1,2}_{q,p}(\Omega_T, x_d^{\alpha p}\omega\, d\mu_2)$ solves \eqref{main-eqn}.

From now on, we fix $k \in \mathbb{N}$, and prove $u_k \in \hat{\sW}^{1,2}_{q,p}(\Omega_T, x_d^{\alpha p}\omega\, d\mu_2)$. Let us denote
\[
D_{R} = (-R^2, \min\{R^2, T\}) \times B_R.\]
Let $R_0>0$ sufficiently large such that the support of $f_k \subset D_{R_0}$.  It follows from \eqref{omega-0}, \eqref{omega-1-0308}, and  H\"{o}lder's inequality that
\[
\|u_k\|_{\hat{\sW}^{1,2}_{q,p}(D_{2R_0}, x_d^{\alpha p} \omega\, d\mu_2)} \leq N(d, p, q, p_1, \alpha, \gamma_0, R_0) \|u_k\|_{\hat{\sW}^{1,2}_{p_1}(D_{2R_0}, x_d^{\alpha p_1}\, d\mu_2)} <\infty.
\]
Then, it remains to prove
\[
\|u_k\|_{\hat{\sW}^{1,2}_{q,p}(\Omega_T\setminus D_{R_0}, x_d^{\alpha p} \omega\, d\mu_2)}  <\infty.
\]
We use a localization technique with the a priori estimate \eqref{main-para.est}. For each $l \in \mathbb{N} \cup \{0\}$, let $\eta_l$ be a smooth function on $(-\infty, T) \times \bR^{d}$ such that
\[
\eta_l \equiv 0 \quad \text{in} \quad D_{2^l R_0}, \quad \eta_l \equiv 1 \quad  (-\infty, T) \times \bR^{d} \setminus D_{2^{l+1} R_0}
\]
and
\[
\|D \eta_l \|_{L_\infty} \leq N_1 2^{-l}, \quad \|\partial_t \eta_l \|_{L_\infty}
+\|D^2\eta_l \|_{L_\infty} \leq N_1 2^{-2l}, \quad \forall \ l \geq 0,
\]
where $N_1$ may depend on $R_0$. Let $w_l = u_k \eta_l$, and we see that $w_l$ solves the equation
\[\mathcal{L} w_l = g \quad \text{in } \Omega_T \quad \text{and} \quad w_l = 0 \quad \text{on } \{x_d =0\},\]
where
\[
g =a_0 u_k \partial_t \eta_l -a_{ij}(D_i u_k D_j \eta_l + D_j u_k D_i \eta_l+u_k D_{ij}\eta_l) - \alpha x_d^{-1} u_k  a_{dj}D_j \eta_l.
\]
Because $w_l \in \hat{\sW}^{1,2}_{p_1}(\Omega_T, x_d^{\alpha p_1} d\mu_2)$, by the estimate \eqref{main-para.est}, it follows that
\[
\begin{split}
\|w_l\|_{\hat{\sW}^{1,2}_{p_1}(\Omega_T, x_d^{\alpha p_1} \omega d\mu_2)} \le N\|g\|_{L_{p_1}(\Omega_T, x_d^{\alpha p_1} \omega d\mu_2)},
\end{split}
\]
where $N$ also depends on $\lambda$. Observe that
\[
\begin{split}
& \|g\|_{L_{p_1}(\Omega_T, x_d^{\alpha p_1} \omega d\mu_2)} \\
& \leq N 2^{-l}\|Du_k\|_{L_{p_1}(D_{2^{l+1}R_0} \setminus D_{2^{l}R_0, x_d^{\alpha p_1} \omega d\mu_2})} + N2^{-2l}\|u_k  \|_{L_{p_1}(D_{2^{l+1}R_0} \setminus D_{2^{l}R_0, x_d^{\alpha p_1} \omega d\mu_2})}\\
&\quad  + N2^{-l}\|\bM^{-1} u_k\|_{L_{p_1}(D_{2^{l+1}R_0} \setminus D_{2^{l}R_0, x_d^{\alpha p_1} \omega d\mu_2})}.
\end{split}
\]
Then, we get
\[
\|u_k\|_{\hat{\sW}^{1,2}_{p_1}(D_{2^{l+2}R_0} \setminus D_{2^{l+1}R_0}, x_d^{\alpha p_1}\omega d\mu_2)} \leq N 2^{-l} \|u_k\|_{\hat{\sW}^{1,2}_{p_1}(D_{2^{l+1}R_0} \setminus D_{2^lR_0}, x_d^{\alpha p_1}\omega d\mu_2)},
\]
where $N>0$ depends also on $R_0$ and $\lambda$, but is independent of $l$. By iterating this estimate, we obtain
\[
\|u_k\|_{\hat{\sW}^{1,2}_{p_1}(D_{2^{l+1}R_0} \setminus D_{2^{l}R_0}, x_d^{\alpha p_1}\omega d\mu_2)} \leq N^l 2^{-l(l-1)/2} \|u_k\|_{\hat{\sW}^{1,2}_{p_1}(D_{2R_0}, x_d^{\alpha p_1}\omega d\mu_2)}, \quad \forall \ l \geq 0.
\]
Finally,  from the inequality above, \eqref{omega-0}, \eqref{omega-1-0308}, and H\"{o}lder's inequality, we obtain
\[
\begin{split}
& \|u_k\|_{\hat{\sW}^{1,2}_{q, p}(\Omega_T \setminus D_{R_0}, x_d^{\alpha p}\omega d\mu_2)}  = \sum_{l=0}^\infty \|u_k\|_{\hat{\sW}^{1,2}_{q, p}(D_{2^{l+1}R_0} \setminus D_{2^{l}R_0}, x_d^{\alpha p}\omega\, d\mu_2)} \\
& \leq  \sum_{l=0}^\infty \|u_k\|_{\hat{\sW}^{1,2}_{p_1}(D_{2^{l+1}R_0} \setminus D_{2^{l}R_0}, x_d^{\alpha p_1} d\mu_2)}     \|\omega_0\|_{L_{1+\varepsilon_1}(\Gamma_{2^{l+1}R_0})}^{1/q} \|\omega_1\|_{L_{1+\varepsilon_1}(B_{2^{l+1}R_0}^+, d\mu_2)}^{1/p}\\
& \leq N \|u_k\|_{\hat{\sW}^{1,2}_{p_1}(D_{R_0}, x_d^{\alpha p_1} d\mu_2)} \sum_{l=0}^\infty N^{l}2^{-l(l-1)/2} N_0^{l(\frac{1}{p} + \frac{1}{q})} <\infty.
\end{split}
\]
The proof is now completed.
\end{proof}
%=====
\begin{proof}[Proof of Theorem \ref{elli-main.theorem}]
Let $\lambda_0$ and $\delta$ be as in Theorem \ref{para-main.theorem}.
It suffices to show the a priori estimate \eqref{main-ell.est} as the existence and uniqueness can be proved in the same way as in the proof of Theorem  \ref{para-main.theorem}.  For a given solution $u \in \hat{\sW}^2_{p}(\bR^{d}_+, x_d^{\alpha p} \omega d\mu_2)$ of \eqref{elli-main-eqn}, let $v(t, x) = \xi(t/n) u(x)$, where $\xi \in C_{0}^\infty((0,1))$. Then, we see that $v \in \hat{\sW}^{1,2}_{p}(\bR^{d+1}_+, x_d^{\alpha p}\omega d\mu_2)$ is a solution of the parabolic equation
\[
v_t - \sL v = g \quad \text{in} \quad \bR^{d+1}_+ \quad \text{with} \quad v =0 \quad \text{on} \quad \{x_d =0\},
\]
where $\sL$ is defined in \eqref{sL} and
\[g(t,x) = \xi(t/n) f(x) + \xi'(t/n) u(x)/n.
\]
By the assumptions in Theorem \ref{elli-main.theorem}, we see that all conditions in Theorem  \ref{para-main.theorem} are satisfied. Then, applying the estimate \eqref{main-para.est} of Theorem  \ref{para-main.theorem} for $v$, and then taking the limit as $n\rightarrow \infty$, we obtain  \eqref{main-ell.est}. See, for example, the proof of \cite[Theorem 1.2]{Dong-Phan-RMI} for details. The theorem is proved.
 \end{proof}

Finally, we give the proof of Corollary \ref{cor2.3}.
\begin{proof}[Proof of Corollary \ref{cor2.3}]
For $k=1,2,\ldots$, we denote  $I_k=(-1+2^{-k},1-2^{-k})$,
$$
Q^k= I_{2k}\times (I_k)^d
\quad \text{and} \quad
Q^k_+= Q^k\cap \Omega_0.
$$
We take a sequence of cutoff functions $\eta_k= \phi_{2k}(t)\prod_{j=1}^d\phi_k(x_j),k=1,2,\ldots,$ where $\phi_k$ satisfies
$$
 \phi_k=1\quad\text{in}\,\,I_{k},\quad \phi_k=0\quad\text{outside}\,\,I_{k+1},
\quad |\phi'_k|\le N2^{k},\quad |\phi''_k|\le N2^{2k}.
$$
Recall the constant $\lambda_0$ from Theorem \ref{para-main.theorem}.
Then it is easily seen that $u\eta_k$ satisfies
\begin{equation} \label{eq2.42}
\left\{
\begin{array}{cccl}
\mathcal{L} (u\eta_k) +\lambda_k c u\eta_k & = & f_k & \quad \text{in} \quad \Omega_{0},    \\
u\eta_k & = & 0 & \quad \text{on} \quad (-\infty, 0) \times \partial \bR^d_+,
\end{array} \right.
\end{equation}
where $\lambda_k\ge \lambda_0\rho_0^{-2}$ is a constant to be specified, $\Omega_0 = (-\infty, 0) \times \bR^d_+$, and
\[
\begin{split}
f_k& =f\eta_k+\lambda_k c u\eta_k+a_0u\eta_t-(a_{ij}+a_{ji})D_iuD_j\eta_k \\
& \quad -a_{ij}uD_{ij}\eta_k- \alpha x_d^{-1} a_{dj} uD_j\eta_k.
\end{split}
\]
It follows from Theorem \ref{para-main.theorem} applied to \eqref{eq2.42} that
\begin{equation} \label{eq2.55}
\begin{split}
A_k&\leq  N\norm{f_k}_{L_{q,p}(\Omega_0, x_d^{p\alpha}\omega\, d\mu_2)}\\
&\le N\|f\|_{L_{q,p}(Q^{k+1}_+, x_d^{p\alpha}\omega\, d\mu_2)}
+N(\lambda_k+2^{2k})\|u\|_{L_{q,p}(Q_+^{k+1}, x_d^{p\alpha}\omega\, d\mu_2)}\\
&\quad + N2^{k}\|\bM^{-1}u\|_{L_{q,p}(Q_+^{k+1}, x_d^{p\alpha}\omega\, d\mu_2)}
+ N2^k\|Du\|_{L_{q,p}(Q^{k+1}_+, x_d^{p\alpha}\omega\, d\mu_2)},
\end{split}
\end{equation}
where
\begin{align*}
A_k&:=
\big\|(u\eta_k)_t|+|DD_{x'}(u\eta_k)| +|\bM^{-1} D_{x'}(u\eta_k)|
+\sqrt{\lambda_k}|D(u\eta_k)|
\\
&\quad  +\sqrt{\lambda_k}|\bM^{-1} u\eta_k|\big\|_{L_{q,p}(\Omega_0, x_d^{p\alpha}\omega\, d\mu_2)} + \norm{D_d(\bM^\alpha D_d (u\eta_k))}_{L_{q,p}(\Omega_0, \omega d\mu_2)},
\end{align*}
and we used the definition of $f_k$ in the last inequality.
From \eqref{eq2.55} and the properties of $\eta_k$, we get
\begin{align}
                    \label{eq4.52}
A_k&\le N2^{k}\lambda_{k+1}^{-1/2}A_{k+1}+N\|f\|_{L_{q,p}(Q^{k+1}_+, x_d^{p\alpha}\omega\, d\mu_2)}\notag\\
&\quad +N(\lambda_k+2^{2k})\|u\|_{L_{q,p}(Q_+^{k+1}, x_d^{p\alpha}\omega\, d\mu_2)}.
\end{align}
We take $\lambda_k=\lambda_0\rho_0^{-2}+(5N2^k)^2$ so that $N2^{k}\lambda_{k+1}^{-1/2}\le 1/5$.
Multiplying both sides of \eqref{eq4.52} by $5^{-k}$ and taking the sum in $k=1,2,\ldots$, we obtain
\begin{align}
                    \label{eq4.53}
\sum_{k=1}^\infty 5^{-k}A_k&\le \sum_{k=1}^\infty 5^{-k-1}A_{k+1}
+N\|f\|_{L_{q,p}(Q^+_1, x_d^{p\alpha}\omega\, d\mu_2)}\nonumber\\
&\quad +N\sum_{k=1}^\infty 5^{-k}(\lambda_k+2^{2k})\|u\|_{L_{q,p}(Q_1^+, x_d^{p\alpha}\omega\, d\mu_2)}.
\end{align}
Note that the summations above are all finite.
By absorbing the first summation on the right-hand side of \eqref{eq4.53} to the left-hand side, we reach
$$
A_1\le N\|f\|_{L_{q,p}(Q^+_1, x_d^{p\alpha}\omega\, d\mu_2)}
+N\|u\|_{L_{q,p}(Q_1^+, x_d^{p\alpha}\omega\, d\mu_2)},
$$
which implies \eqref{eq2.17}. The corollary is proved.
\end{proof}

We conclude the paper with the following remark.
\begin{remark} \label{final-remark} It is possible to study the class of \eqref{main-eqn} and \eqref{elli-main-eqn} with the additional zeroth order terms of the form $b u/x_d^2$. For example, let us consider
we consider the equation
\begin{equation} \label{main-eqn-bis}
\left\{
\begin{array}{cccl}
\mathcal{L}v(t,x)  - \frac{b}{x_d^2} v(t,x) & = &  f(t, x)  &\quad \text{in} \quad  \Omega_T\\
 v & =& 0 & \quad \text{on}\quad (-\infty, T) \times \partial \bR^d_+,
\end{array} \right.
\end{equation}
where $b$  is a constant and $\mathcal{L}$ is defined in \eqref{L.def}. We also assume
\begin{equation} \label{fina-remark-cond}
a_{dd}\equiv 1 \quad \text{and} \quad a_{dj} = a_{jd} \quad \text{for all} \quad  j=1,2,\ldots, d-1.
\end{equation}
Following \cite{MNS21}, we define $u(t,x) = x_d^{\beta} v(t,x)$. Then formally, $u$ satisfies
\[
\hat{\mathcal{L}} u := a_0 u_t - a_{ij} D_{ij} u - \frac{\alpha-2\beta}{x_d}a _{dj} D_{j} u + \lambda c u- \frac{\beta^2 -(\alpha -1) \beta + b}{x_d^2}u=x_d^{\beta}f
. \]
We shall choose a $\beta$ such that
\begin{equation} \label{beta-eqn}
\beta^2 -(\alpha -1) \beta + b =0.
\end{equation}
When $b < (\alpha-1)^2/4$, \eqref{beta-eqn} has two real roots
\[
\beta_1=\frac{\alpha -1 + \sqrt{(\alpha-1)^2 -4b}}{2} \quad \text{and} \quad \beta_2=\frac{\alpha -1 - \sqrt{(\alpha-1)^2 -4b}}{2}.
\]
Denote $\alpha_{i} = \alpha -2\beta_{i}$ for $i =1, 2$, so that
$$
\alpha_1  = 1 - \sqrt{(\alpha-1)^2 -4b} <1 \quad \text{and} \quad \alpha_{2} = 1 + \sqrt{(\alpha-1)^2 -4b} >1. $$
Now we take $\beta=\beta_1$, $\mu_2(x_d) = x_d^{\gamma_0}$, and $\gamma_0 \in (-1, \alpha_1 -1]$.
By \eqref{fina-remark-cond} and Remark \ref{rem2.1}, we see that $u\in \hat{\sW}^{1,2}_{q,p}(\Omega_T, x_d^{\alpha_1 p} \omega d\mu_2)$ is a solution to
\begin{equation} \label{main-eqn-new}
\left\{
\begin{array}{cccl}
 \hat{\mathcal{L}}u(t,x)  & =  & x_d^{\beta}f(t, x) &\quad \text{in} \quad  \Omega_T, \\
 u & = & 0 & \quad \text{on} \quad (-\infty, T) \times \partial \bR^d_+
 \end{array}\right.
\end{equation}
if and only if $v$ is a strong solution to \eqref{main-eqn-bis} and it satisfies
\begin{align*}
&\|v\|+\|\bM^{-1}v\|+\|v_t\|+\|Dv\|+\|\bM^{-1}D_{x'}v\|+\|DD_{x'}v\|
\\
&\quad +\|D_d^2 v+\alpha \bM^{-1} D_dv+b\bM^{-2} v\|<\infty,
\end{align*}
where $\|\cdot\|=\|\cdot\|_{L_{q,p}(\Omega_T, x_d^{(\alpha_1+\beta_1) p} \omega d\mu_2)}$.
We then apply Theorem \ref{para-main.theorem} to obtain the unique solvability of \eqref{main-eqn-new} and the estimate for $u \in \hat{\sW}^{1,2}_{q,p}(\Omega_T, x_d^{\alpha_1 p} \omega d\mu_2)$. Then,  by changing back to $v$, we can derive the corresponding result for \eqref{main-eqn-bis}.
\end{remark}
%=====
\section*{Acknowledgement}

After we finished  an earlier version of this manuscript and posted it online, we were informed by Prof. Metafune about his preprint \cite{MNS21}, where they obtained the maximal regularity estimate for \eqref{eq1.52} and \eqref{eq5.51} with an additional zeroth-order term $bu/x_d^2$, where $b$ is a constant, and  $p\in (\alpha p-1,2p-1)$ when $b=0$. Their proof is completely different from ours. We thank Prof. Metafune for letting us know this reference.

\end{document}